\documentclass[11pt,reqno]{amsart}
\usepackage{amsmath,amssymb,amsthm,amsrefs}
\usepackage[linktocpage=true,colorlinks,citecolor=magenta,linkcolor=blue,urlcolor=magenta]{hyperref}
\usepackage[margin=1.1in]{geometry}
\newtheorem{thm}{Theorem}[section]
\newtheorem{pro}{Proposition}[section]
\newtheorem{lemma}{Lemma}[section]
\newtheorem{corollary}{Corollary}[section]
\theoremstyle{remark}
\newtheorem{remark}{Remark}[section]

\numberwithin{equation}{section}
\allowdisplaybreaks
\hypersetup{pdftitle={Alexandrov-Fenchel inequalities for convex bodies in hyperbolic space},
  pdfauthor={Yong Wei and Xinxiu Zhang},
  pdfsubject={Sharp comparisons between hyperbolic quermassintegrals of geodesically convex bodies},
  pdfkeywords={quermassintegrals, Alexandrov-Fenchel inequalities, hyperbolic space, support measures}}

\begin{document}
\lefthyphenmin=2
\righthyphenmin=3
\title[Alexandrov--Fenchel inequalities in hyperbolic space]
  {Alexandrov--Fenchel inequalities for convex bodies in hyperbolic space}

\author[Y. Wei]{Yong Wei}
\author[X. Zhang]{Xinxiu Zhang}
\address{School of Mathematical Sciences, University of Science and Technology of China, Hefei 230026, P.R. China}
\email{\href{mailto:yongwei@ustc.edu.cn}{yongwei@ustc.edu.cn}}
\email{\href{mailto:zhangxinxiu@mail.ustc.edu.cn}{zhangxinxiu@mail.ustc.edu.cn}}

\date{\today}
\subjclass[2020]{52A39, 52A55, 53C42}
\keywords{Alexandrov--Fenchel inequality, quermassintegral, convex body, hyperbolic space, support measure}
\thanks{This work was supported by the National Key Research and Development Program of China (2021YFA1001800), National Natural Science Foundation of China (No.~12531002) and the Fundamental Research Funds for the Central Universities.}

\begin{abstract}
Let $\Omega$ be a nonempty compact geodesically convex set in $\mathbb H^{n+1}$. For $0\leq\ell<k\leq n$, we prove that $W_k(\Omega)\geq W_k(B)$, where $B$ is a geodesic ball of possibly zero radius with $W_\ell(B)=W_\ell(\Omega)$. If $B$ has positive radius, equality holds precisely when $\Omega$ is a geodesic ball. At radius zero, equality holds precisely when $\dim\Omega\leq n-k$. The proof reduces the comparison of consecutive quermassintegrals to a variational problem. For a minimizer with interior, one-sided variations give an integral inequality that we localize to the polar faces over exposed points. We show that a curvature quotient is concave on the conormal cones and use continuous approximations on the polar faces to apply this inequality to the quotient. Newton's inequalities for the curvature measures then give the required upper bound for the total Gauss curvature. For a lower-dimensional minimizer we obtain the bound by restriction to its totally geodesic span.
\end{abstract}

\maketitle
\tableofcontents

\section{Introduction}\label{sec:intro}

Among Euclidean bodies of given volume, the isoperimetric inequality says that a ball has the least surface area. For convex bodies the  Alexandrov--Fenchel
inequalities give a corresponding comparison for every pair of
quermassintegrals. At fixed quermassintegral $W_\ell$, a ball minimizes $W_k$ for each $k>\ell$. We refer to Schneider~\cite{schneider}*{\S\S7.3--7.4}
for more details on the Euclidean Alexandrov--Fenchel inequalities.

In a space form, the quermassintegral $W_j$ of a convex body is,
up to normalization, the invariant measure of the totally geodesic
$j$-planes meeting the body. See Santal\'o~\cite{santalo}*{Ch.~17}
and Solanes~\cite{solanes}. Gao, Hug and
Schneider~\cite{gaohugschneider}*{p.~165} asked which
integral-geometric functionals attain an extremum precisely at
balls among spherical or hyperbolic convex bodies of given
volume. Their two-point symmetrization gives the fixed-volume
comparison for the measure of the hyperplanes meeting a body.
The Euclidean inequalities suggest a further question. If a lower
quermassintegral $W_\ell$ is prescribed, does a geodesic ball
minimize every $W_k$ with $k>\ell$?

In the sphere $\mathbb S^{n+1}$, Chen, Guan, Li and
Scheuer~\cite{cgls}*{Conjecture~1.1} conjectured these inequalities
for convex bodies. Luo, Wei and Zhou~\cite{lwzsphere}*{Theorem~1.1}
proved the full family, including the equality case, for smooth
convex bodies contained in an open hemisphere. For earlier
inequalities and curvature flows in the sphere, see
\cites{makowskischeuer,weixiong,chensun,crscheuer,panscheuer}.

In hyperbolic space $\mathbb H^{n+1}$ we use the normalization of Andrews, Chen
and Wei~\cite{acw} on the quermassintegrals, in which $W_0(\Omega)=|\Omega|$ and
$nW_1(\Omega)=|\partial\Omega|$ for convex bodies with nonempty
interior. Wang and Xia~\cite{wangxia} proved all the comparisons
for smooth horospherically convex bodies, whose principal
curvatures satisfy $\kappa_i\geq1$. Along their flow, $W_\ell$
is constant, $W_k$ is nonincreasing, and the hypersurfaces
converge to a geodesic sphere. Related quermassintegral-preserving
flows were studied by Andrews and Wei~\cite{aw} and by Hu, Li
and Wei~\cite{hlw}. The second author~\cite{zhangh}
extended the inequalities and their equality characterization to
possibly nonsmooth horospherically convex bodies.

Under weaker assumptions on the boundary, the known comparisons in hyperbolic space depend on the pair $(\ell,k)$. Li, Wei and Xiong~\cite{lwx}
proved the $(1,3)$ comparison for smooth star-shaped two-convex
hypersurfaces. Brendle, Guan and Li~\cite{bgl}
obtained the $(1,2)$ comparison for smooth star-shaped
hypersurfaces of positive mean curvature and the $(\ell,n)$
comparisons for smooth strictly convex bodies. Both results are
stated in Guan and Li~\cite{guanlisurvey}*{Theorem~4.2,
parts~(1) and~(2)}. Pan and Yang~\cite{panyang}*{Theorem~1.8}
proved the $(1,2)$ comparison for smooth outward-minimizing
boundaries when $2\leq n\leq6$. Hu and Li~\cite{huli}*{Theorem~1.3}
proved the $(\ell,n-1)$ comparisons for smooth strictly convex
bodies. For $(0,n)$, Wei, Yang and Zhou~\cite{wyz}*{Corollary~1.3}
gave a curvature flow proof for smooth weakly convex bodies.

For a smooth convex hypersurface of dimension $n\geq2$, the
Gauss equation gives $K_{ij}=-1+\kappa_i\kappa_j$ in each
principal two-plane. Nonnegative intrinsic sectional curvature
thus lies between horospherical convexity and geodesic convexity.
For odd $k\geq3$, Ge, Wang and Wu~\cite{gewangwu} compared $W_k$
with $W_1$ under horospherical convexity, and Hu and Li~\cite{hlodd}
obtained the same comparisons under nonnegative sectional
curvature. Andrews, Chen and Wei~\cite{acw}*{Corollary~1.5}
proved all comparisons with volume under positive sectional
curvature. For $n\geq3$, Luo, Wei, Zhang and
Zhou~\cite{luo}*{Theorem~1.1} proved the comparisons under
nonnegative sectional curvature whenever $k+\ell\leq n$.
For $1\leq\ell<k\leq n$ with $k+\ell>n$, their curvature
quotient flow can develop negative sectional curvature from
initial data of positive sectional curvature while the
hypersurface remains smooth and strictly convex~\cite{luo}*{Theorem~1.4}.

In this paper, we prove the full family of comparisons under geodesic convexity alone.

\begin{thm}\label{thm:main}
Let $n\geq1$ and let $\Omega\subset\mathbb H^{n+1}$ be a nonempty
compact geodesically convex set. Write $f_j(r)=W_j(B(r))$ for a
geodesic ball $B(r)$ of radius $r>0$, and set $f_j(0)=0$.
For every $0\leq\ell<k\leq n$,
\begin{equation}\label{eq:main}
 W_k(\Omega)\geq f_k\bigl(f_\ell^{-1}(W_\ell(\Omega))\bigr).
\end{equation}
If $W_\ell(\Omega)>0$, equality holds if and only if $\Omega$ is
a geodesic ball. If $W_\ell(\Omega)=0$, equality holds if and only if
$\dim\Omega\leq n-k$, where $\dim\Omega$ is the dimension of its
totally geodesic span.
\end{thm}

In particular, among convex bodies of prescribed positive volume,
the geodesic ball uniquely minimizes each $W_k$, $1\leq k\leq n$,
up to isometry. This answers the fixed-volume question of Gao, Hug
and Schneider~\cite{gaohugschneider} for the hyperbolic Crofton functionals.

\subsection{Outline of the proof}

The inequalities \eqref{eq:main} say that the comparison radii
$r_j(K)=f_j^{-1}(W_j(K))$ are nondecreasing in $j$. By Hausdorff
approximation and continuity of the quermassintegrals, it suffices
to prove the inequalities for bodies with nonempty interior. Suppose that
$\ell$ is the least index for which $r_{\ell+1}(K)<r_\ell(K)$
for some convex body $K$. We minimize the perturbed deficit
\begin{equation*}
 J(K)=W_{\ell+1}(K)-\psi(W_\ell(K))+\varepsilon W_n(K),
\end{equation*}
where $\psi$ is a bounded truncation of
$f_{\ell+1}\circ f_\ell^{-1}$. Writing $\varphi=f_{\ell+1}\circ f_\ell^{-1}$, we set $\psi(t)=M\bigl(1-e^{-\varphi(t)/M}\bigr)$.
For $M$ sufficiently large and $\varepsilon>0$ sufficiently small,
we have $J(K)<0$ at the chosen counterexample $K$.
Since $W_n$ controls the diameter, compactness up to isometries
gives a minimizer $\Omega_*$. Put $r=r_\ell(\Omega_*)>0$ and
$c=\coth r$. The choice of $\ell$ and the recurrence for the total
normalized mean curvatures $T_j$ give
$T_\ell(\Omega_*)<T_\ell(B(r))$. By \eqref{eq:top}, the comparison
with $W_{n-1}$ due to Hu and Li~\cite{huli} gives
$T_n(\Omega_*)\geq T_n(B(r))$ when $\ell<n-1$. For $\ell=n-1$
these two quantities agree by the choice of $r$. Hence
\begin{equation}\label{s1.Tn}
 T_n(\Omega_*)\geq T_n(B(r))
 =c^{n-\ell}T_\ell(B(r))
 >c^{n-\ell}T_\ell(\Omega_*).
\end{equation}

The truncation gives $0<\gamma<c$, where
$\gamma=\psi'(W_\ell(\Omega_*))$. If $\Omega_*$ were smooth with
positive principal curvatures, the first variation and Newton's
inequalities for the normalized mean curvatures $E_j$ would give
\begin{equation*}
 E_{\ell+1}+\varepsilon E_n=\gamma E_\ell,
 \qquad
 \frac{E_n}{E_\ell}
 \leq\left(\frac{E_{\ell+1}}{E_\ell}\right)^{n-\ell}
 \leq\gamma^{n-\ell}.
\end{equation*}
Integration would contradict the strict lower bound \eqref{s1.Tn} for $T_n(\Omega_*)$.

Assume now that $\Omega_*$ has nonempty interior. Its boundary may be
nonsmooth, so we work with curvature measures. In their study of a
constrained curvature flow in Euclidean space, Andrews and
the first author~\cite{awjdg}*{proof of Lemma~6.3} showed that two curvature
measures of a limiting convex body are proportional and used
Schneider's theorem to conclude that the limit is a ball. Here we
need a local inequality between curvature measures of different orders
on $\partial\Omega_*$. At a nonsmooth
point, projection to the boundary combines the contributions from
different supporting normals. Write $V_{n-j}/\binom nj=a_j\lambda$,
where $\lambda$ is the sum of these normalized curvature measures.
In Euclidean space the local Steiner formula expresses the
curvature measures in terms of mixed Hessian measures of the
distance function. We apply the mixed Hessian inequality of Dinew
and Lu~\cite{dinew}*{Theorem~3.9} in the real form of
Lemma~\ref{lem:real-mixed}. Comparing the Klein metric with
constant metrics on small balls then gives the hyperbolic
inequalities $a_j^2\geq a_{j-1}a_{j+1}$. We also show that the
indices $j$ for which $a_j(x)>0$ form an interval at
$\lambda$-almost every $x$. It therefore suffices to prove
$a_\ell>0$ and $a_{\ell+1}\leq\gamma a_\ell$ on $\{a_n>0\}$.

We use one-sided convex variations in Klein coordinates of the form
$\Omega_t=\bigcup_{x\in\Omega_*}(x+tL_x)$ to obtain an integral
inequality for products $\chi(x)b(y)$. Here $\chi$ is nonnegative
and convex, and $b$ is concave polyhedral and nonnegative on the
polar body. Convex functions that vanish on $\Omega_*$ outside
small neighborhoods of exposed points localize the inequality.
Disintegrating the support measures over the contact point and
approximating on each polar face gives the inequality for every
continuous nonnegative concave function on that face, outside one
$\lambda$-null set of contact points.

Write $\mu_{j,x}$ for the fiber measure of mass $a_j(x)$ on the
unit conormal section $S_x$. Using disintegration on the faces of
a convex function \cites{caravenna,ciosmakF}, we show that for
$\lambda$-almost every $x$ with $a_n(x)>0$, the point $x$ is exposed
and $\mu_{n,x}$ is equivalent to spherical measure on $S_x$.
The constant function $1$ on the polar face gives $a_\ell(x)>0$.
To estimate $a_{\ell+1}(x)$, we express the densities
$d\mu_{j,x}/d\mu_{n,x}$ in terms of the Hessian of the support
function. On the relative interior of the conormal cone these
Hessians have a common image wherever they exist, and their
inverses on this image are superadditive. It follows that
$d\mu_{\ell+1,x}=Q\,d\mu_{\ell,x}$ on $S_x$, where $Q$ is
nonnegative, concave and homogeneous of degree one. The quotient
$Q$ is continuous on the relative interior of the cone. We set $Q=0$
on its relative boundary. Approximating the restriction of $Q$
to the polar face by continuous nonnegative concave functions
gives the inequality for $Q$, and homogeneity transfers it to
$S_x$. By Cauchy--Schwarz,
\begin{equation*}
 a_{\ell+1}^2
 \leq a_\ell\int_{S_x}Q^2\,d\mu_{\ell,x}
 =a_\ell\int_{S_x}Q\,d\mu_{\ell+1,x}
 \leq\gamma a_\ell a_{\ell+1}.
\end{equation*}
The successive ratios $a_{j+1}/a_j$ are nonincreasing for
$\ell\leq j<n$ on $\{a_n>0\}$, so $a_n\leq\gamma^{n-\ell}a_\ell$.
Since $T_j=\int a_j\,d\lambda$, integration gives
$T_n(\Omega_*)\leq\gamma^{n-\ell}T_\ell(\Omega_*)$.

If $\Omega_*$ has empty interior, we restrict $J$ to its totally
geodesic span, where $\Omega_*$ has nonempty interior. Applying the
interior estimate there and using the restriction formulas in
Lemma~\ref{lem:restrict} gives the same upper bound in the ambient
space. In either case we obtain a contradiction. This proves
\eqref{eq:main}.

To characterize equality of \eqref{eq:main} at positive radius, we show that $r_m=r_{m+1}>0$
implies $r_{m-1}=r_m$ for $1\leq m<n$. Suppose instead that
$r_{m-1}<r_m=r_{m+1}$. We perturb the deficit for the pair $(m,m+1)$
by $\delta(W_{m-1}-h(W_m)+\eta W_n)$, where $\delta,\eta>0$ and
$h$ is a bounded truncation of $f_{m-1}\circ f_m^{-1}$.
The perturbed minimizers retain a positive lower bound for
$f_{m-1}(r_m)-W_{m-1}$ as their $(m,m+1)$ deficit tends to zero.
By the recurrence \eqref{eq:rec}, $T_m$ then stays strictly below
$T_m(B(r_m))$ by a fixed amount. On the other hand, the curvature
estimate gives $T_m(B(r_m))-T_m\leq C\delta$ for small $\delta$.
This proves the implication. Iteration gives $r_0=r_1$,
so rigidity follows from the hyperbolic isoperimetric inequality.
At radius zero the equality cases follow directly from the Crofton formula.

\subsection{Organization of the paper}

Section~\ref{sec:prelim} recalls the formulas for quermassintegrals
and support measures and proves a local inequality for mixed
Hessian measures. In Section~\ref{sec:reduction} we reduce the
problem to minimizing a perturbed deficit.
Section~\ref{sec:newton} establishes Newton's inequalities for
the curvature measures. In Section~\ref{sec:variation} we derive
an integral inequality on the polar faces over exposed points.
In Section~\ref{sec:fibre} we apply this inequality to the
curvature quotient by approximation and bound the total Gauss
curvature of a minimizer with interior. Section~\ref{sec:proof} treats
lower-dimensional minimizers and characterizes equality for all compact
geodesically convex sets.

\subsection{AI Declaration}  
The authors used OpenAI's ChatGPT (6 Astra and 6.2 Astra) to assist with exploring proof strategies, searching the literature, developing and checking mathematical arguments and calculations, and improving the exposition. The authors take full responsibility for all mathematical statements, proofs, references, and conclusions in the final manuscript.

\section{Preliminaries}\label{sec:prelim}

\subsection{Quermassintegrals in hyperbolic space}
Let $\Omega$ be a nonempty compact convex set in $\mathbb H^{n+1}$. Its quermassintegrals are defined using the invariant measure $dL_k$ on the space $\mathcal L_k$ of totally geodesic $k$-planes (see~\cite{solanesthesis}*{Definition~2.2.1} and~\cite{solanes}). Write $\omega_j$ for the area of the unit sphere $\mathbb S^j$. In the normalization of~\cite{acw}*{(1.8)}, $W_0(\Omega)=|\Omega|$ is the hyperbolic volume and
\begin{equation*}
 W_k(\Omega)=\frac{\omega_{k-1}\cdots\omega_0}{\omega_{n-1}\cdots\omega_{n-k}}
       \int_{\mathcal L_k}\chi(L_k\cap\Omega)\,dL_k,
       \qquad 1\le k\le n,
\end{equation*}
where $\chi$ is the Euler characteristic. In the convention of~\cites{gallego,wangxia}, the quermassintegral of order $k$ is $\frac{n+1-k}{n+1}W_k(\Omega)$.  For convex $\Omega$, the integrand equals one if $L_k$ meets $\Omega$ and zero otherwise. The quermassintegrals are monotone under inclusion and continuous with respect to the Hausdorff distance \cites{gallego,acw}.

For a smooth closed hypersurface $M=\partial\Omega$ with outer unit normal $\nu$, let $\kappa_1,\ldots,\kappa_n$ be the eigenvalues of the Weingarten map $v\mapsto\overline\nabla_v\nu$. The normalized $k$-th mean curvature is $E_k(\kappa)=\binom nk^{-1}\sigma_k(\kappa)$, where $\sigma_k$ is the $k$-th elementary symmetric function and $E_0=1$. For a smooth family $M_t=\partial\Omega_t$ of convex hypersurfaces with outward normal speed $v$ and area element $d\mu_t$, the first variation formula (see~\cite{wangxia}*{Proposition~3.1} and~\cite{acw}*{Lemma~2.3}) gives
\begin{equation}\label{eq:smoothvar}
 \frac{d}{dt}W_k(\Omega_t)=\int_{M_t}vE_k(\kappa)\,d\mu_t,
 \qquad 0\le k\le n.
\end{equation}
Write $f_k(r)=W_k(B(r))$ for the geodesic ball of radius $r$. For $r>0$, its boundary has area $\omega_n\sinh^n r$ and all principal curvatures equal to $\coth r$. Taking $v=1$ in \eqref{eq:smoothvar} gives
\begin{equation}\label{eq:ball}
 f_k'(r)=\omega_n\sinh^{n-k}r\,\cosh^k r,\qquad f_k(0)=0,
 \qquad 0\le k\le n.
\end{equation}

Working in the Klein model, one readily sees that every nonempty compact
convex set in $\mathbb H^{n+1}$ can be approximated in the Hausdorff metric
by smooth strictly convex bodies. By the continuity of
quermassintegrals,  the following lemma is an immediate
consequence of \cite{huli}*{Theorem~1.3}.
\begin{lemma}\label{lem:upperindex}
Let $n\ge2$ and let $\Omega\subset\mathbb H^{n+1}$ be a nonempty compact
convex set. Then
\begin{align*}
 W_{n-1}(\Omega)
 &\ge f_{n-1}\circ f_\ell^{-1}\bigl(W_\ell(\Omega)\bigr),
 \qquad 0\le\ell\le n-2.
\end{align*}
\end{lemma}

\subsection{Support measures and curvature measures}
For a nonempty compact geodesically convex set $\Omega$, let
$f_\Omega$ be the metric projection from $\mathbb H^{n+1}$ onto 
$\Omega$. For $z\notin\Omega$, let $\nu_z$ be the unit tangent vector
at $f_\Omega(z)$ pointing toward $z$, and set
$F_\Omega(z)=(f_\Omega(z),\nu_z)\in T^1\mathbb H^{n+1}$.
The unit normal bundle of $\Omega$ is
$\mathcal N_\Omega=F_\Omega(\mathbb H^{n+1}\setminus\Omega)$.
For $\varepsilon>0$, let $\Omega^\varepsilon=\{z:d(z,\Omega)\leq\varepsilon\}$
be the outer parallel set at distance $\varepsilon$. Kohlmann
\cite{kohlmann} proved that there exist Radon measures
$\Phi_k(\Omega,\cdot)$ on $T^1\mathbb H^{n+1}$ satisfying the local
Steiner formula
\begin{equation}\label{eq:tube}
 \operatorname{Vol}\bigl(F_\Omega^{-1}(\beta)\cap(\Omega^\varepsilon\setminus\Omega)\bigr)
 =\sum_{k=0}^n\Phi_k(\Omega,\beta)
              \int_0^\varepsilon \sinh^{n-k}s\,\cosh^k s\,ds
\end{equation}
for Borel sets $\beta\subset T^1\mathbb H^{n+1}$ and $\varepsilon>0$. For any $k\in\{0,\dots,n\}$, we call $\Phi_k(\Omega,\cdot)$ the $k$-th support measure of $\Omega$. Let $\Pi:T^1\mathbb H^{n+1}\to\mathbb H^{n+1}$ be the bundle projection. Each support measure $\Phi_k(\Omega,\cdot)$ is concentrated on $\mathcal N_\Omega$. The curvature measures are the projections
$V_k(\Omega,\cdot)=\Pi_\#\Phi_k(\Omega,\cdot)$.
If $\partial\Omega$ is smooth, then for every Borel set
$\beta\subset\mathbb H^{n+1}$,
\begin{align*}
&V_k(\Omega, \beta)=\binom{n}{k}\int_{\beta\cap\partial\Omega}{E}_{n-k} d\mu,  \quad 0\leq k \leq n.
\end{align*}
 For the representation of curvature measures by differential forms, see~\cite{valuation}*{\S2.3}. 

For convex bodies with nonempty interior, the curvature measures
are weakly continuous under Hausdorff convergence
\cite{kohlmann1994}*{Theorem~2.4}. See also
\cite{veronelli}*{\S2}. For the support measures, the proof of
\cite{zhangh}*{Theorem~2.4} also applies when the limit has empty
interior. Indeed, if $K_i\to K$ and $z\notin K$, uniqueness of
the nearest point gives $F_{K_i}(z)\to F_K(z)$. When $K$ has
empty interior, both $K$ and $\partial K^s$ have zero ambient
volume for every $s>0$. Dominated convergence therefore gives
convergence of the local tube integrals against continuous
compactly supported test functions. The coefficient functions
of $s$ in \eqref{eq:tube} are linearly independent, so the
support measures converge weakly at each order.


We write
\begin{equation*}
 T_j(\Omega)=\binom nj^{-1}\Phi_{n-j}(\Omega,\mathcal N_\Omega),
 \qquad 0\le j\le n,
\end{equation*}
so that $T_j(\Omega)=\int_{\partial\Omega}E_j(\kappa)\,d\mu$ when $\partial\Omega$ is smooth. The total mean curvatures satisfy
\begin{align}
 T_j(\Omega)&=(n-j)W_{j+1}(\Omega)+jW_{j-1}(\Omega)
                    &&(0\le j<n),\label{eq:rec}\\
 T_n(\Omega)&=\omega_n+nW_{n-1}(\Omega).\label{eq:top}
\end{align}
At $j=0$ the second term in \eqref{eq:rec} is omitted. For smooth convex bodies these are the identities in~\cite{gallego}*{(5)}, written in the normalization of~\cite{acw}*{(1.10)--(1.12)}. Smooth approximation and the continuity of both sides extend them to every nonempty compact convex set.
\subsection{Monge--Amp\`ere and mixed Hessian measures}
\label{subsec:ma}

For a finite convex function $v$ on $\mathbb R^{n+1}$, its
Alexandrov Monge--Amp\`ere measure is
\begin{equation*}
 \operatorname{MA}(v)(B)=|\partial v(B)|,
 \qquad \partial v(B)=\bigcup_{x\in B}\partial v(x),
\end{equation*}
where $B$ is Borel, $\partial v(x)$ is the subdifferential,
and $|\cdot|$ denotes Lebesgue measure. This is a nonnegative
Radon measure, with $\operatorname{MA}(v)=\det D^2v\,dx$ for
$v\in C^2$ \cite{clm}*{Lemma~4.1 and Theorem~4.2}.
We write $\operatorname{MA}(v_1,\ldots,v_{n+1})$ for its symmetric
multilinear polarization, normalized by
$\operatorname{MA}(v,\ldots,v)=\operatorname{MA}(v)$.
These mixed measures are nonnegative and jointly weakly
continuous under local uniform convergence
\cite{clm}*{Theorem~4.3}.
We write $v[k]$ for $k$ repetitions of $v$ and set
$q(x)=|x|^2/2$, so that $\operatorname{MA}(q)=dx$.

For finite convex functions $v_1,\ldots,v_m$ and
$1\le m\le n+1$, define
\begin{equation*}
 H_m(v_1,\ldots,v_m)
 =\operatorname{MA}(v_1,\ldots,v_m,q[n+1-m]),
 \qquad H_m(v)=H_m(v[m]).
\end{equation*}
For $v\in C^2$, this gives
$H_m(v)=\binom{n+1}{m}^{-1}\sigma_m(D^2v)\,dx$,
where $\sigma_m$ is the $m$th elementary symmetric polynomial
of the eigenvalues; see \cite{clm}*{\S4.2}.

The following lemma is a consequence of the local mixed
Hessian inequality of Dinew and Lu
\cite{dinew}*{Theorem~3.9}. We include the reduction to real
convex functions.

\begin{lemma}\label{lem:real-mixed}
Let $2\le m\le n+1$ and let $u_1,\ldots,u_m$ be finite convex
functions on $\mathbb R^{n+1}$. Suppose that $\mu$ is a nonnegative
Radon measure and $f_i\ge0$ are Borel functions such that
$H_m(u_i)\ge f_i\mu$ for $1\le i\le m$. Then
\begin{equation}\label{newton:real-mixed}
 H_m(u_1,\ldots,u_m)
 \ge\left(\prod_{i=1}^m f_i\right)^{1/m}\mu.
\end{equation}
\end{lemma}

\begin{proof}
For a finite convex function $v$, the function
$\widetilde v(x+\mathrm i y)=v(x)$ is locally bounded and
plurisubharmonic. Set $\omega=dd^c\widetilde q$. Let $dy$ denote Lebesgue measure in the imaginary variables and
define $b_{n+1}>0$ by $\omega^{n+1}=b_{n+1}\,dx\,dy$.
For smooth convex $v$, the complex Hessian matrices of
$\widetilde v$ and $\widetilde q$ are $\tfrac14 D^2v$ and
$\tfrac14 I$, respectively. Hence
\begin{equation*}
 (dd^c\widetilde v)^{n+1}=(\det D^2v)\,\omega^{n+1}.
\end{equation*}
Polarizing this identity with $v_1,\ldots,v_m$ and
$n+1-m$ copies of $q$ gives
\begin{equation}\label{prelim:real-complex}
 dd^c\widetilde v_1\wedge\cdots\wedge dd^c\widetilde v_m
 \wedge\omega^{n+1-m}
 =b_{n+1}H_m(v_1,\ldots,v_m)\otimes dy .
\end{equation}

For general finite convex $v_i$, take smooth convex
approximations decreasing to $v_i$ and converging locally
uniformly. The left side of \eqref{prelim:real-complex}
converges weakly by \cite{bt}*{Theorem~2.1}, and the right 
side converges weakly by \cite{clm}*{Theorem~4.3}. Thus
\eqref{prelim:real-complex} holds as an equality of Radon
measures on $\mathbb C^{n+1}$.

Let $E=\{\prod_i f_i>0\}$ and
$\nu=b_{n+1}(\mu|_E)\otimes dy$. The assumptions give
\begin{equation*}
 (dd^c\widetilde u_i)^m\wedge\omega^{n+1-m}
 \ge f_i(x)\,\nu,\qquad 1\le i\le m.
\end{equation*}
Since $f_1>0$ on $E$, we have
$\nu\ll(dd^c\widetilde u_1)^m\wedge\omega^{n+1-m}$.
Since $\widetilde u_1$ is locally bounded, its complex Hessian
measure $(dd^c\widetilde u_1)^m\wedge\omega^{n+1-m}$ vanishes
on $m$-polar sets \cite{dinew}*{\S2}. By absolute continuity,
$\nu$ also vanishes on $m$-polar sets.
On every relatively compact complex ball, the functions
$\widetilde u_i$ are bounded and $m$-subharmonic.
Thus \cite{dinew}*{Theorem~3.9} gives locally
\begin{equation*}
 dd^c\widetilde u_1\wedge\cdots\wedge dd^c\widetilde u_m
 \wedge\omega^{n+1-m}
 \ge\left(\prod_{i=1}^m f_i(x)\right)^{1/m}\nu.
\end{equation*}
Using \eqref{prelim:real-complex} and testing on product sets
gives \eqref{newton:real-mixed} with $\mu|_E$ in place of $\mu$.
The right side vanishes outside $E$, proving the claim.
\end{proof}
\subsection{Relation between Euclidean curvature measures and mixed Hessian measures}
\label{subsec:euclidean-curvature}

We recall the relation between Euclidean curvature measures
and the Hessian measures of the distance function
\cite{hug}*{(1.3)}.
Let $\Omega\subset\mathbb R^{n+1}$ be a nonempty compact convex set,
and put $d(x)=\operatorname{dist}(x,\Omega)$.
Write $V_k^{\mathrm{eucl}}=V_k^{\mathrm{eucl}}(\Omega,\cdot)$,
$0\le k\le n$, for its Euclidean curvature measures, normalized
so that
\begin{equation*}
 V_k^{\mathrm{eucl}}=\binom nk E_{n-k}\,d\mu
\end{equation*}
when $\partial\Omega$ is smooth, with the curvatures and area
element computed in the Euclidean metric.

For $x\in\partial\Omega$, the subdifferential $\partial d(x)$
is the intersection of the Euclidean outer normal cone
$\operatorname{Nor}(\Omega,x)$ with the closed unit ball.
Thus, for every Borel set $B\subset\partial\Omega$ and $t>0$,
\begin{equation*}
 \partial(q+td)(B)
 =\{x+v:x\in B,\ v\in\operatorname{Nor}(\Omega,x),\ |v|\le t\}.
\end{equation*}
The right-hand side is the union of $B$ and the set of points
outside $\Omega$ at distance at most $t$ whose nearest point
in $\Omega$ lies in $B$.
Since $B$ has zero Lebesgue measure, the local Steiner formula
\cite{schneider}*{\S4.2}, in the above normalization, gives
\begin{equation*}
 \operatorname{MA}(q+td)(B)
 =\sum_{j=1}^{n+1}\frac{t^j}{j}
 V_{n+1-j}^{\mathrm{eucl}}(B).
\end{equation*}
Comparing coefficients of $t^j$ in the mixed Monge--Amp\`ere
expansion yields
\begin{equation}\label{newton:steiner}
 \operatorname{MA}(d[j],q[n+1-j])\big|_{\partial\Omega}
 =\frac{V_{n+1-j}^{\mathrm{eucl}}}
 {(n+1)\binom n{j-1}},
 \qquad 1\le j\le n+1.
\end{equation}

\section{Reduction to a variational problem}\label{sec:reduction}

The inequalities in~\eqref{eq:main} are equivalent to
\begin{equation*}
 f_0^{-1}(W_0(\Omega))\le\cdots\le f_n^{-1}(W_n(\Omega)).
\end{equation*}
Suppose that one of these comparisons fails. We truncate the corresponding ball comparison function and add $\varepsilon W_n$ to the deficit. The truncation gives the derivative bound in~\eqref{eq:gamma}. The term $\varepsilon W_n$ controls the diameter of a minimizing sequence and contributes a positive multiple of the Gauss curvature measure to the first variation.

\begin{pro}\label{prop:min}
 Suppose that $\ell$ is the least index in $\{0,\ldots,n-1\}$ for which some compact geodesically convex body $\Omega_0$ with nonempty interior in $\mathbb H^{n+1}$ satisfies $W_{\ell+1}(\Omega_0)<\varphi(W_\ell(\Omega_0))$, where $\varphi=f_{\ell+1}\circ f_\ell^{-1}$. There exist $M>0$ and $\varepsilon>0$ such that the functional
\begin{equation}\label{eq:penalty}
 J(\Omega)=W_{\ell+1}(\Omega)-M\bigl(1-\exp[-\varphi(W_\ell(\Omega))/M]\bigr)
                       +\varepsilon W_n(\Omega)
\end{equation}
has a negative minimum among nonempty compact geodesically convex sets. If $\Omega_*$ is a minimizer and $W_\ell(B(r))=W_\ell(\Omega_*)$, then $r>0$ and
\begin{equation}\label{eq:strict-lower}
 T_n(\Omega_*)\ge T_n(B(r))
       =(\coth r)^{n-\ell}T_\ell(B(r))
       >(\coth r)^{n-\ell}T_\ell(\Omega_*).
\end{equation}
The derivative of $t\mapsto M(1-\exp[-\varphi(t)/M])$ at $t=W_\ell(\Omega_*)$ is
\begin{equation}\label{eq:gamma}
 0<\gamma=e^{-f_{\ell+1}(r)/M}\coth r<\coth r.
\end{equation}
\end{pro}
\begin{proof}
Choose a convex body $K$ for which the comparison at $\ell$ fails. As $M\to\infty$ and $\varepsilon\downarrow0$,
\begin{equation*}
 J(K)\longrightarrow W_{\ell+1}(K)-\varphi(W_\ell(K))<0.
\end{equation*}
Fix $M,\varepsilon>0$ for which $J(K)<0$.

Let $I$ be a diameter segment of $\Omega$. The Crofton formula~\cite{solanesthesis}*{Theorem~2.2.2} gives $W_n(I)=(\omega_{n-1}/n)\operatorname{length}(I)$. Since $I\subset\Omega$, monotonicity gives
\begin{equation*}
 J(\Omega)\ge\varepsilon W_n(\Omega)-M
 \ge\varepsilon W_n(I)-M
 =\varepsilon\frac{\omega_{n-1}}{n}\operatorname{diam}(\Omega)-M.
\end{equation*}
Thus $J$ is bounded below, and a minimizing sequence with $J<0$
has uniformly bounded diameters. Since $J$ is invariant under
hyperbolic isometries, we may assume that each set in the sequence
contains the origin. In the Klein model, these are Euclidean convex
sets contained in a common compact subball. Blaschke's selection
theorem~\cite{schneider}*{Theorem~1.8.7} therefore yields a subsequence
converging in the Euclidean, and hence hyperbolic, Hausdorff metric
to a nonempty compact convex set. By continuity of the
quermassintegrals, $J$ attains its infimum at this limit, and
$\inf J<0$.

Let $\Omega_*$ be a minimizer and  $r=f_\ell^{-1}(W_\ell(\Omega_*))$. Since $J(\Omega_*)<0$, we have
\begin{equation*}
 0\le W_{\ell+1}(\Omega_*)+\varepsilon W_n(\Omega_*)
 < M\bigl(1-e^{-f_{\ell+1}(r)/M}\bigr)
 \le f_{\ell+1}(r).
\end{equation*}
In particular, $r>0$ and $W_{\ell+1}(\Omega_*)<f_{\ell+1}(r)$. If $\ell>0$, the comparison at $\ell-1$ holds for all convex bodies with nonempty interior by the choice of $\ell$. Approximation as in the proof of Lemma~\ref{lem:upperindex} extends it to $\Omega_*$. Thus
\begin{equation*}
 f_{\ell-1}^{-1}(W_{\ell-1}(\Omega_*))
 \le f_\ell^{-1}(W_\ell(\Omega_*))=r.
\end{equation*}
It follows from~\eqref{eq:rec} that
\begin{align*}
 T_\ell(B(r))-T_\ell(\Omega_*)
 &=(n-\ell)\bigl(f_{\ell+1}(r)-W_{\ell+1}(\Omega_*)\bigr)\\
 &\quad+\ell\bigl(f_{\ell-1}(r)-W_{\ell-1}(\Omega_*)\bigr)>0,
\end{align*}
where the second term is omitted when $\ell=0$.

If $\ell<n-1$, Lemma~\ref{lem:upperindex} gives
$W_{n-1}(\Omega_*)\geq f_{n-1}(r)$. If $\ell=n-1$, then
$W_{n-1}(\Omega_*)=f_{n-1}(r)$ by the choice of $r$. Thus~\eqref{eq:top} gives
\begin{equation*}
 T_n(\Omega_*)-T_n(B(r))
 =n\bigl(W_{n-1}(\Omega_*)-f_{n-1}(r)\bigr)\ge0.
\end{equation*}
Since $T_j(B(r))=\omega_n\sinh^n r\,(\coth r)^j$, these two comparisons give~\eqref{eq:strict-lower}.

Finally, by~\eqref{eq:ball},
\begin{equation*}
 0<\gamma
 =e^{-f_{\ell+1}(r)/M}\frac{f_{\ell+1}'(r)}{f_\ell'(r)}
 =e^{-f_{\ell+1}(r)/M}\coth r<\coth r. \qedhere
\end{equation*}
\end{proof}

In view of~\eqref{eq:strict-lower} and~\eqref{eq:gamma},
a contradiction will follow if we prove
\begin{equation*}
 T_n(\Omega_*)\le\gamma^{n-\ell}T_\ell(\Omega_*).
\end{equation*}
For a minimizer with nonempty interior, this estimate is established
in Proposition~\ref{prop:total}. For a minimizer with empty interior,
we apply that proposition in its totally geodesic span, using the
restriction formulas in Section~\ref{sec:dimension}.
We first establish Newton's inequalities for the curvature measures.

\section{Newton's inequalities for curvature measures}\label{sec:newton}

For a nonsmooth convex body, the curvature measures may have singular
parts with respect to surface area measure. We prove Newton's inequalities
for their Radon--Nikodym derivatives with respect to a common dominating
measure.

\begin{pro}\label{prop:newton}
Let $n\ge1$ and let $\Omega$ be a nonempty compact
convex set in $\mathbb R^{n+1}$ or $\mathbb H^{n+1}$. 
Write $V_k=V_k(\Omega,\cdot)$ for its curvature measures
in the respective ambient space, normalized so that
$V_{n-j}=\binom nj E_j\,d\mu$ when $\partial\Omega$
is smooth, and set
\begin{equation*}
 \lambda=\sum_{i=0}^n\frac{V_{n-i}}{\binom ni},
 \qquad
 a_j=\frac{1}{\binom nj}\frac{dV_{n-j}}{d\lambda},
 \quad 0\le j\le n.
\end{equation*}
Then, for $\lambda$-almost every $x\in\partial\Omega$,
\begin{equation}\label{newton:projected}
 \begin{aligned}
 a_j(x)^2&\ge a_{j-1}(x)a_{j+1}(x)
 &&(1\le j<n),\\
 a_i(x)a_k(x)>0&\quad\Longrightarrow\quad a_j(x)>0
 &&(0\le i<j<k\le n).
 \end{aligned}
\end{equation}
\end{pro}

\subsection{The Euclidean case}
\label{subsec:newton-euclidean}

We first prove the Euclidean case of
Proposition~\ref{prop:newton}.
\begin{proof}[Euclidean case of Proposition~\ref{prop:newton}]
Let $\Omega\subset\mathbb R^{n+1}$ be a nonempty compact
convex set, put $d(x)=\operatorname{dist}(x,\Omega)$,
and set
\begin{equation*}
 M_j=\operatorname{MA}(d[j],q[n+1-j])\big|_{\partial\Omega},
 \qquad 0\le j\le n+1.
\end{equation*}
Since $\operatorname{MA}(q)=dx$ and $\partial\Omega$ has zero
Lebesgue measure, $M_0=0$.
Set $\Lambda=\sum_{j=0}^{n+1}M_j$, and choose nonnegative
Borel representatives $h_j=dM_j/d\Lambda$ with $h_0=0$
and $\sum_{j=0}^{n+1}h_j=1$ $\Lambda$-almost everywhere.
For $2\le m\le n+1$, set
\begin{equation*}
 P_m(x,t)=\sum_{j=0}^m\binom mj h_j(x)t^{m-j},
 \qquad x\in\partial\Omega,\quad t>0.
\end{equation*}

Fix $2\le m\le n+1$ and $t>0$.
By the definitions and multilinearity,
$H_m(d)\ge h_m\Lambda$ and
$H_m(d+tq)\ge P_m(\cdot,t)\Lambda$,
where the densities are extended by zero off $\partial\Omega$.
Applying Lemma~\ref{lem:real-mixed} with $\mu=\Lambda$
to $m-1$ copies of $d$ and one copy of $d+tq$,
and restricting the resulting measure inequality to
$\partial\Omega$, we obtain by multilinearity
\begin{equation}\label{newton:polynomial}
 \bigl(h_m(x)+t h_{m-1}(x)\bigr)^m
 \ge h_m(x)^{m-1}P_m(x,t)
\end{equation}
for $\Lambda$-almost every $x\in\partial\Omega$.

Take $t_k=1/k$, $k\ge1$.
Then, for $\Lambda$-almost every $x\in\partial\Omega$,
\eqref{newton:polynomial} holds with $t=t_k$ for every
$k\ge1$ and $2\le m\le n+1$, and
$\sum_{j=0}^{n+1}h_j(x)=1$.
If $h_m(x)>0$ at such a point $x$, expanding the numerator
in powers of $t_k$ gives
\begin{equation*}
 \begin{aligned}
 0\le{}&
 \frac{h_m(x)\left(1+t_k\frac{h_{m-1}(x)}{h_m(x)}\right)^m
       -P_m(x,t_k)}
      {\binom m2t_k^2}\\
 &\longrightarrow
 \frac{h_{m-1}(x)^2}{h_m(x)}-h_{m-2}(x).
 \end{aligned}
\end{equation*}
Thus
\begin{equation*}
 h_{m-1}(x)^2\ge h_m(x)h_{m-2}(x),
 \qquad 2\le m\le n+1,
\end{equation*}
which also holds when $h_m(x)=0$.

If $h_m(x)>0$ and $h_{m-1}(x)=0$ for some
$2\le m\le n+1$, taking $t=t_1=1$ in
\eqref{newton:polynomial} gives $P_m(x,1)\le h_m(x)$,
forcing $h_0(x)=\cdots=h_{m-1}(x)=0$.
If a zero occurred between two positive entries, the first
positive entry after that zero would contradict this implication.
Consequently,
\begin{equation*}
 h_i(x)h_k(x)>0\quad\Longrightarrow\quad h_j(x)>0,
 \qquad 1\le i<j<k\le n+1.
\end{equation*}
The set of positive indices is nonempty since
$\sum_{j=0}^{n+1}h_j(x)=1$.

Set
\begin{equation*}
 \lambda_{\mathrm{eucl}}
 =\sum_{j=0}^n\frac{V_{n-j}^{\mathrm{eucl}}}{\binom nj}.
\end{equation*}
By \eqref{newton:steiner},
\begin{equation*}
 \Lambda=\frac{\lambda_{\mathrm{eucl}}}{n+1},\qquad \quad 
 h_{j+1}(x)
 =\frac{1}{\binom nj}
   \frac{dV_{n-j}^{\mathrm{eucl}}}
        {d\lambda_{\mathrm{eucl}}}(x),
 \qquad 0\le j\le n,
\end{equation*}
where the density identities hold
$\lambda_{\mathrm{eucl}}$-almost everywhere on $\partial\Omega$.
This proves \eqref{newton:projected} in Euclidean space.
\end{proof}

\subsection{The hyperbolic case}

To prove Proposition~\ref{prop:newton} in hyperbolic space, we compare
the support measures of a convex set $\Omega$ in the Klein ball $B^{n+1}$
for the hyperbolic metric $g$ and for a constant positive definite
metric $G$. In Klein coordinates, we write $p(v)$ for the
evaluation of a covector $p\in(\mathbb R^{n+1})^*$ on a vector
$v\in\mathbb R^{n+1}$. The Euclidean support function is
$h_\Omega:(\mathbb R^{n+1})^*\to\mathbb R$, given by
$h_\Omega(p)=\max_{z\in\Omega}p(z)$. For $x\in\partial\Omega$,
the outer conormal cone is
\begin{equation}\label{eq:conormal-cone}
 \mathcal N_x
 =\bigl\{p\in(\mathbb R^{n+1})^*:p(z-x)\leq0
                   \text{ for all }z\in\Omega\bigr\}
 =\bigl\{p\in(\mathbb R^{n+1})^*:h_\Omega(p)=p(x)\bigr\}.
\end{equation}
For $p\ne0$, write $[p]=\{tp:t>0\}$. The initial tangent to the
geodesic from $x$ to $z\ne x$ is a positive multiple of $z-x$ in
Klein coordinates. The first variation and convexity of squared
distance along geodesics therefore identify the nonzero rays in
\eqref{eq:conormal-cone} with the rays of the metric duals of
outer unit normals for both $g$ and $G$.

Consider the cosphere bundle \cite{valuation}*{\S2.1}
\begin{equation*}
 S^*B^{n+1}=(T^*B^{n+1}\setminus\{0\})/\mathbb R_{>0}.
\end{equation*}
For $m\in\{g,G\}$, write $|p|_{m_x^{-1}}=\sqrt{m_x^{-1}(p,p)}$ for
the dual norm. The map
\begin{equation*}
 (x,[p])\longmapsto
 \left(x,\frac{m_x^{-1}p}{|p|_{m_x^{-1}}}\right)
\end{equation*}
identifies $S^*B^{n+1}$ with the $m$-unit tangent bundle and sends
outer conormal rays to the corresponding outer unit normals. Here
$p/|p|_{m_x^{-1}}$ is the unit conormal for the dual metric. We use
the same notation for the support measures transported to $S^*B^{n+1}$
and write $\Pi(x,[p])=x$ for the base projection.

Write $\overline\nabla$ for the Levi-Civita connection of $g$ and
$D$ for the Euclidean connection.
In Klein coordinates (see Besau and Werner~\cite{besauwerner}),
\begin{equation}\label{eq:connection}
 \begin{aligned}
 g_x(v,w)
 &=\frac{v\cdot w}{1-|x|^2}
   +\frac{(x\cdot v)(x\cdot w)}{(1-|x|^2)^2},\\
 \overline\nabla_v w-D_v w
 &=\frac{(x\cdot v)w+(x\cdot w)v}{1-|x|^2}.
 \end{aligned}
\end{equation}
For the constant metric $G$ on $\mathbb R^{n+1}$, write $\Phi_k^G$
for the support measures and $V_k^G=\Pi_\#\Phi_k^G$ for their
projections. We normalize $\Phi_k^G$ by \eqref{eq:tube}, using
$G$-distance and $G$-volume and replacing $\sinh s,\cosh s$ by
$s,1$, respectively. Thus in the same bundle representation
$\Phi_k^G=\omega_{n-k}\Lambda_k^G$, where $\Lambda_k^G$ is the
support measure in Schneider's normalization \cite{schneider}*{\S4.2}. Viewed as measures on $S^*B^{n+1}$ 
via the above homeomorphisms, $\Phi_k$ and $\Phi_k^G$ remain weakly
continuous under Hausdorff convergence of nonempty compact convex sets
(Section~\ref{sec:prelim} and \cite{schneider}*{\S4.2}, respectively).

\begin{lemma}\label{lem:metric}
Let $\Omega$ be a nonempty compact convex set in the Klein ball and
let $U_0\Subset B^{n+1}$ be open. There exist $c,\delta_0>0$,
depending only on $n$ and bounds for $g$ and $g^{-1}$ on $U_0$,
with the following property. If $U\subset U_0$ is open and $G$ is a
constant positive definite metric with
$\sup_U\|g-G\|\le\delta\le\delta_0$ in the Euclidean operator norm,
then
\begin{equation}\label{newton:freeze}
 e^{-c\delta}\Phi^G_{n-j}\le\Phi_{n-j}\le e^{c\delta}\Phi^G_{n-j}
 \quad\text{on }\Pi^{-1}(U),\qquad 0\le j\le n.
\end{equation}
In particular, the same bounds hold for $V_{n-j}$ and $V^G_{n-j}$
on $U$.
\end{lemma}

\begin{proof}
Suppose first that $\Omega$ has smooth boundary and nonempty
interior. Fix $z\in\partial\Omega\cap U$ and a nonzero outer
conormal $p$ at $z$. For tangent vector fields $v,w$,
the difference $\overline\nabla_v w-D_v w$ in
\eqref{eq:connection} is tangent to $\partial\Omega$ and
hence annihilated by $p$. Since $D$ is also the
Levi-Civita connection of $G$, we have
\begin{equation*}
 \mathrm{II}^g(v,w)
 =-\frac{p(D_v w)}{|p|_{g_z^{-1}}}
 =\frac{|p|_{G^{-1}}}{|p|_{g_z^{-1}}}
  \mathrm{II}^G(v,w).
\end{equation*}

The assumption on $G$ and the uniform bounds for $g$
and $g^{-1}$ imply
\begin{equation*}
 e^{-\eta}G\le g\le e^\eta G
 \quad\text{on }U,
 \qquad \eta=C\delta,
\end{equation*}
where $C$ depends only on the metric bounds in the statement.
Consequently,
\begin{equation*}
 e^{-\eta/2}
 \le\frac{|p|_{G^{-1}}}{|p|_{g_z^{-1}}}
 \le e^{\eta/2}.
\end{equation*}

For every nonzero $v\in T_z\partial\Omega$, we have
\begin{equation*}
 \frac{\mathrm{II}^g(v,v)}{g_z(v,v)}
 =\frac{|p|_{G^{-1}}}{|p|_{g_z^{-1}}}
  \frac{G(v,v)}{g_z(v,v)}
  \frac{\mathrm{II}^G(v,v)}{G(v,v)}.
\end{equation*}
Since $\mathrm{II}^G$ is positive semidefinite,
the preceding metric comparisons and the min--max
principle yield
\begin{equation*}
 e^{-3\eta/2}\kappa_i^G
 \le\kappa_i^g
 \le e^{3\eta/2}\kappa_i^G,
 \qquad 1\le i\le n,
\end{equation*}
where the principal curvatures are arranged in
nondecreasing order.

Moreover,
\begin{equation*}
 e^{-n\eta/2}d\mu_G
 \le d\mu_g
 \le e^{n\eta/2}d\mu_G.
\end{equation*}
Since all principal curvatures are nonnegative, it follows
that
\begin{equation*}
 \begin{aligned}
 e^{-(n+3j)\eta/2}
 \binom nj E_j^G\,d\mu_G
 &\le \binom nj E_j^g\,d\mu_g~\le e^{(n+3j)\eta/2}
 \binom nj E_j^G\,d\mu_G,
 \end{aligned}
\end{equation*}
for $0\le j\le n$. Identifying the smooth support measures
through the common conormal rays and using
$(n+3j)/2\le2n$ proves \eqref{newton:freeze}
with $c=2nC$.

For a general nonempty convex subset $\Omega$, choose smooth strictly convex bodies
$\Omega_m$ converging to $\Omega$ in the Hausdorff metric
within a common compact subball. The preceding estimates
hold uniformly in $m$. By weak continuity of the support
measures, they pass to the limit when tested against
nonnegative functions in $C_c(\Pi^{-1}(U))$.
This proves \eqref{newton:freeze} for $\Omega$.
The corresponding inequalities for curvature measures
follow by projection under $\Pi$.
\end{proof}

For every compact set $K\Subset B^{n+1}$ in the Klein model,
the same argument with fixed comparison constants shows that
the support measures of each order associated with $g$ and
any constant positive definite metric $G$ are mutually
absolutely continuous on $\Pi^{-1}(K)$.
The same holds for any two constant positive definite
metrics and for their projected curvature measures on $K$.

\begin{proof}[Hyperbolic case of Proposition~\ref{prop:newton}]
Let $\Omega\subset B^{n+1}$ be a nonempty compact convex set,
and choose an open set $U_0\Subset B^{n+1}$ containing $\Omega$.
Fix a sequence $\delta_k\downarrow0$ with
$0<\delta_k<\delta_0$. By the continuity of $g$ and the
compactness of $\partial\Omega$, for each $k$ there is a finite
open cover by open balls 
\begin{equation*}
 \partial\Omega\subset\bigcup_{\alpha=1}^{N_k}U_{k,\alpha},
 \qquad U_{k,\alpha}\subset U_0,
\end{equation*}
such that
\begin{equation*}
 \sup_{U_{k,\alpha}}\|g-G_{k,\alpha}\|
 \le\delta_k,
 \qquad
 G_{k,\alpha}=g_{x_{k,\alpha}},
\end{equation*}
where $x_{k,\alpha}$ is the center of $U_{k,\alpha}$.

Fix $k$ and $\alpha$, and write $U=U_{k,\alpha}$ and
$G=G_{k,\alpha}$. Set
\begin{equation*}
 \lambda_G=\sum_{i=0}^n\frac{V^G_{n-i}}{\binom ni},
 \qquad
 b_j^G=\frac{1}{\binom nj}
       \frac{dV^G_{n-j}}{d\lambda_G},
 \quad 0\le j\le n.
\end{equation*}
Since $G$ is constant and positive definite, the linear map
$x\mapsto G^{1/2}x$ is an isometry from
$(\mathbb R^{n+1},G)$ to standard Euclidean space.
Thus the Euclidean case of Proposition~\ref{prop:newton}
applies to the densities $b_j^G$.

By the projected measure comparison in
\eqref{newton:freeze}, $\lambda_G$ and $\lambda$ are
equivalent on $U$. Hence
\begin{equation*}
 a_j^G=\frac{1}{\binom nj}\frac{dV^G_{n-j}}{d\lambda}
      =\frac{d\lambda_G}{d\lambda}\,b_j^G,
 \qquad 0\le j\le n,
\end{equation*}
$\lambda$-almost everywhere on $U$.
Since the factor $d\lambda_G/d\lambda$ is positive and
independent of $j$, both assertions of
\eqref{newton:projected} hold for $a_j^G$ on $U$.

Moreover, \eqref{newton:freeze} gives
\begin{equation*}
 e^{-c\delta_k}a_j^G
 \le a_j\le e^{c\delta_k}a_j^G,
 \qquad 0\le j\le n,
\end{equation*}
for $\lambda$-almost every $x\in U$. Consequently,
\begin{equation*}
 a_{j-1}a_{j+1}
 \le e^{2c\delta_k}a_{j-1}^Ga_{j+1}^G~\le e^{2c\delta_k}(a_j^G)^2
 \le e^{4c\delta_k}a_j^2,
 \qquad 1\le j<n.
\end{equation*}
Since the cover is finite, this inequality holds for
$\lambda$-almost every $x\in\partial\Omega$, simultaneously
for all $1\le j<n$. Taking the intersection of these
full-measure sets over $k$, we may let $k\to\infty$ to obtain
\begin{equation*}
 a_j^2\ge a_{j-1}a_{j+1},
 \qquad 1\le j<n,
\end{equation*}
for $\lambda$-almost every $x\in\partial\Omega$.

For the second assertion, fix one of the finite covers
above. On each $U_{k,\alpha}$, the two-sided comparison
implies
\begin{equation*}
 a_j(x)>0\quad\Longleftrightarrow\quad
 a_j^{G_{k,\alpha}}(x)>0,
 \qquad 0\le j\le n,
\end{equation*}
for $\lambda$-almost every $x\in U_{k,\alpha}$.
The Euclidean case therefore shows that the positive
indices of $(a_0(x),\ldots,a_n(x))$ have no gaps for
$\lambda$-almost every $x\in\partial\Omega$.
This completes the proof.
\end{proof}

Let $\Omega_*$ be a minimizer as in Proposition~\ref{prop:min}.
For $\lambda$-almost every $x$ with $a_\ell(x)a_n(x)>0$,
\eqref{newton:projected} gives positivity of the intermediate
densities and shows that $a_{j+1}(x)/a_j(x)$ is nonincreasing
for $\ell\le j<n$. Hence
\begin{equation*}
 \frac{a_n}{a_\ell}
 =\prod_{j=\ell}^{n-1}\frac{a_{j+1}}{a_j}
 \le\left(\frac{a_{\ell+1}}{a_\ell}\right)^{n-\ell}.
\end{equation*}
Thus, to establish
$T_n(\Omega_*)\le\gamma^{n-\ell}T_\ell(\Omega_*)$,
it suffices to prove that $a_\ell>0$ and
$a_{\ell+1}\le\gamma a_\ell$ $\lambda$-almost everywhere
on $\{a_n>0\}$.

\section{Convex variations and an inequality on polar faces}
\label{sec:variation}

We derive first variation formulas for nonempty compact convex sets
from the local tube formula. For local minimizers with nonempty
interior, we use the polar parametrization of the normal bundle
to obtain integral inequalities and localize them to the polar
faces over exposed points.

\subsection{First variations of compact convex sets}
For the curvature form identities underlying first variation formulas
in space forms, see~\cite{fu1998}*{Lemma~3.1}. We use the local
Steiner formula in a form that also applies to the one-sided
variations in Proposition~\ref{prop:variation}. For a related
distance-function argument for continuous Steiner symmetrization in
Euclidean space, see~\cite{ghomi2026}*{Appendix~B, pp.~12--13}.

\begin{lemma}\label{lem:ambient}
Let $\Omega\subset\mathbb H^{n+1}$ be a nonempty compact convex set,
and let $\varphi_t$ be a smooth family of diffeomorphisms of a
neighborhood of $\Omega$ onto their images, with $\varphi_0=\mathrm{id}$.
Suppose that $\varphi_t(\Omega)$ is geodesically convex for $|t|$
sufficiently small. Then, for $0\le j\le n$,
\begin{equation}\label{eq:ambientvar}
 \left.\frac{d}{dt}\right|_0 W_j(\varphi_t(\Omega))
 =\frac{1}{\binom nj}\int_{\mathcal N_\Omega}
    \langle\dot\varphi_0(x),\nu\rangle_g\,\Phi_{n-j}(\Omega,d(x,\nu)),
\end{equation}
where $g$ is the hyperbolic metric and $\nu$ is the outer unit
normal at the incident pair $(x,\nu)$.
\end{lemma}

\begin{proof}
Let $c_j(s)=\sinh^j s\,\cosh^{n-j}s$. For every nonempty compact
convex set $K$ and every $s>0$, \eqref{eq:tube}, \eqref{eq:rec}
and \eqref{eq:top} give
\begin{align*}
 \frac{d}{ds}\operatorname{Vol}(K^s)
 &=\sum_{j=0}^n\binom nj c_j(s)T_j(K)\\
 &=\omega_n c_n(s)+\sum_{j=0}^n\binom nj c_j'(s)W_j(K),
\end{align*}
where $c_j'=j c_{j-1}+(n-j)c_{j+1}$, with zero terms omitted.
Since $W_0(K)=\operatorname{Vol}(K)$, integration gives
\begin{equation*}
 \operatorname{Vol}(K^s)=\operatorname{Vol}(B(s))
       +\sum_{j=0}^n\binom nj c_j(s)W_j(K).
\end{equation*}

Let $K_t$ be a family of nonempty compact convex sets with
$K_0=\Omega$ and $d_H(K_t,\Omega)=O(t)$ as $t\downarrow0$. Set
$d_t(z)=\min_{a\in K_t}d_g(z,a)$ and suppose that, for some bounded
Borel function $v$ on $\mathcal N_\Omega$,
\begin{equation*}
 \left.\frac{d}{dt}\right|_{0+}d_t(z)
 =-v(F_\Omega(z))
 \qquad (z\notin\Omega).
\end{equation*}
For $\rho\in C_c^\infty((0,\infty))$, Fubini's theorem gives
\begin{equation*}
 I_\rho(t):=\int_0^\infty \rho(s)\operatorname{Vol}(K_t^s)\,ds
 =\int_{\mathbb H^{n+1}}\int_{d_t(z)}^\infty
           \rho(s)\,ds\,d\operatorname{vol}(z).
\end{equation*}
The difference quotients of the inner integral are bounded in
absolute value by $\|\rho\|_\infty d_H(K_t,\Omega)/t$ and vanish
outside a fixed compact set. For $z\in\Omega$, they also vanish
for all sufficiently small $t$, since $\rho$ vanishes near zero.
By the local Steiner formula \eqref{eq:tube},
the pushforward of the volume measure on
$\mathbb H^{n+1}\setminus\Omega$ under the map
$z\mapsto(F_\Omega(z),d_0(z))$ is
\begin{equation*}
 \sum_{j=0}^n
 \Phi_{n-j}(\Omega,\cdot)\otimes(c_j(s)\,ds).
\end{equation*}
Hence, by dominated convergence and Fubini's theorem,
\begin{align*}
 I_\rho'(0+)
 &=\int_{\mathbb H^{n+1}\setminus\Omega}
       \rho(d_0(z))v(F_\Omega(z))\,d\operatorname{vol}(z)\\
 &=\sum_{j=0}^n
       \left(\int_0^\infty\rho(s)c_j(s)\,ds\right)
       \int_{\mathcal N_\Omega}v\,d\Phi_{n-j}(\Omega).
\end{align*}
Choose a nonzero $\rho_0\in C_c^\infty((0,\infty))$ with
$\rho_0\ge0$ and take $\rho=\rho_0c_i$ for $0\le i\le n$.
The matrix
\begin{equation*}
 \left(\int_0^\infty\rho_0(s)c_i(s)c_j(s)\,ds\right)_{i,j=0}^n
\end{equation*}
is positive definite, since $c_j(s)/\cosh^n s=\tanh^j s$.
Substituting the tube formula into the difference quotients of
$I_{\rho_0c_i}$ and using the invertibility of the matrix above,
we obtain
\begin{equation*}
 \left.\frac{d}{dt}\right|_{0+}W_j(K_t)
 =\frac{1}{\binom nj}\int_{\mathcal N_\Omega}
       v\,d\Phi_{n-j}(\Omega),
 \qquad 0\le j\le n.
\end{equation*}

Now take $K_t=\varphi_t(\Omega)$. Fix $z\notin\Omega$, write
$(x,\nu)=F_\Omega(z)$, and choose $a_t\in\Omega$ such that
$d_t(z)=d_g(z,\varphi_t(a_t))$. The functions
$a\mapsto d_g(z,\varphi_t(a))$ converge uniformly on $\Omega$ to
$a\mapsto d_g(z,a)$. The latter has the unique minimizer $x$, so
$a_t\to x$. Taylor's formula gives
\begin{equation*}
 d_g(z,\varphi_t(a))=d_g(z,a)+tL(a)+O(t^2),
 \qquad
 L(a)=\left.\frac{\partial}{\partial t}\right|_{t=0}
 d_g(z,\varphi_t(a)),
\end{equation*}
uniformly for $a\in\Omega$. Comparing the values at $a_t$ and $x$
and using $a_t\to x$, we obtain
\begin{equation*}
 \left.\frac{d}{dt}\right|_0 d_t(z)
 =L(x)=-\langle\dot\varphi_0(x),\nu\rangle_g,
\end{equation*}
where the last equality follows from the first variation of distance.
Since $d_H(K_t,\Omega)=O(|t|)$, applying the right derivative
formula to $\varphi_t(\Omega)$ and $\varphi_{-t}(\Omega)$ proves
\eqref{eq:ambientvar}.
\end{proof}

To allow the normal speed to depend on both the contact point and
the conormal, we assign a compact convex set $L_x$ of velocities to
each point $x$ in Klein coordinates. The set update
$\bigcup_{x\in\Omega}(x+tL_x)$ also occurs in set evolution equations
(see~\cite{panasyuk1980}*{(4), p.~431 of the Russian original} and
\cite{dudagehriglorenz}*{Proposition~2.1(2), p.~764}). The support
function of $L_x$, evaluated at the hyperbolic unit conormal, gives
the one-sided normal speed.

\begin{pro}\label{prop:variation}
Let $\Omega$ be a nonempty compact convex set in the Klein model
of $\mathbb H^{n+1}$. Let
$F:\Omega\times(\mathbb R^{n+1})^*\to\mathbb R$ be continuous,
concave in the first variable, and convex and positively homogeneous
of degree one in the second variable. For each $x\in\Omega$, let
$L_x\subset\mathbb R^{n+1}$ be the nonempty compact convex set with
support function $F(x,\cdot)$, and set
\begin{equation*}
 \Omega_t=\{x+tq:x\in\Omega,\ q\in L_x\},\qquad t\geq0.
\end{equation*}
For all sufficiently small $t\geq0$, $\Omega_t$ is a nonempty
compact geodesically convex set in the Klein ball. For $0\leq j\leq n$,
\begin{equation}\label{var:first}
 \left.\frac{d}{dt}\right|_{0+}W_j(\Omega_t)
 =\frac{1}{\binom nj}\int_{\mathcal N_\Omega}
 F(x,p_{x,\nu})\,\Phi_{n-j}(\Omega,d(x,\nu)),
\end{equation}
where $p_{x,\nu}=g_x(\nu,\cdot)$ is the hyperbolic unit conormal
at $(x,\nu)\in\mathcal N_\Omega$, expressed in Klein coordinates.
\end{pro}

\begin{proof}
Let $\mathcal G=\{(a,q):a\in\Omega,\ q\in L_a\}$. The graph
$\mathcal G$ is compact by continuity of $F$ on the product of
$\Omega$ with the Euclidean unit sphere in the dual space, and
convex by concavity of $F$ in the first variable. Its image under
$(a,q)\mapsto a+tq$ is the nonempty compact convex set $\Omega_t$.

For small $t$, both $\Omega$ and $\Omega_t$ lie in a fixed compact
subball of the Klein model, and $d_H(\Omega_t,\Omega)=O(t)$.

Fix $z\notin\Omega$ and let $(x,\nu)=F_\Omega(z)$. Set
$d_t(z)=\min_{a\in\Omega_t}d_g(z,a)$ and $d_z(a)=d_g(z,a)$. Then
\begin{equation*}
 d_t(z)=\min_{(a,q)\in\mathcal G} d_z(a+tq).
\end{equation*}
Since $d_z$ is smooth near $\Omega$ and $\mathcal G$ is compact,
the objective and its $t$-derivative are jointly continuous for
$|t|$ small and $(a,q)\in\mathcal G$. At $t=0$, the set of
minimizers is $\{x\}\times L_x$. By Danskin's theorem
\cite{bonnansshapiro}*{Theorem~4.1} and the first variation of
distance,
\begin{equation}\label{var:distvel}
 \begin{aligned}
 -\left.\frac{d}{dt}\right|_{0+}d_t(z)
 &=-\min_{q\in L_x}Dd_z(x)[q]\\
 &=\max_{q\in L_x}p_{x,\nu}(q)=F(x,p_{x,\nu}).
 \end{aligned}
\end{equation}
The function $(x,\nu)\mapsto F(x,p_{x,\nu})$ is bounded and
continuous on $\mathcal N_\Omega$. Applying the one-sided variation formula from the proof of Lemma~\ref{lem:ambient} with $K_t=\Omega_t$ and
$v(x,\nu)=F(x,p_{x,\nu})$ gives \eqref{var:first}.
\end{proof}

\subsection{Localized first variations}

In this subsection, we assume that $\Omega$ has nonempty interior
and choose the Klein origin in $\operatorname{int}\Omega$.
The Klein coordinates identify $T^*B^{n+1}$ with
$B^{n+1}\times(\mathbb R^{n+1})^*$. For $x\in\partial\Omega$,
the conormal cone in the fiber over $x$ is $\{x\}\times\mathcal N_x$,
with $\mathcal N_x$ as in \eqref{eq:conormal-cone}. For nonzero
$p\in\mathcal N_x$, write $(x,[p])$ for the conormal ray
$\{(x,tp):t>0\}$.

Let $h_\Omega$ denote the Euclidean support function of $\Omega$
in Klein coordinates. The polar body of $\Omega$ is
\begin{equation*}
 \Omega^\circ=\{y\in(\mathbb R^{n+1})^*:h_\Omega(y)\leq1\}.
\end{equation*}
Since $h_\Omega(p)>0$ for $p\ne0$, each conormal ray $(x,[p])$
has a unique representative $(x,y)$ with $h_\Omega(y)=1$,
given by $y=p/h_\Omega(p)$.
By \eqref{eq:conormal-cone}, $y(x)=h_\Omega(y)=1$. This gives the
incidence set
\begin{equation*}
 \mathcal I=\{(x,y)\in\Omega\times\Omega^\circ:y(x)=1\}.
\end{equation*}
The corresponding parametrization of the hyperbolic unit normal
bundle is
\begin{equation*}
 J_\Omega:\mathcal N_\Omega\longrightarrow\mathcal I,\qquad
 J_\Omega(x,\nu)
 =\left(x,\frac{p_{x,\nu}}{h_\Omega(p_{x,\nu})}\right),
 \qquad p_{x,\nu}=g_x(\nu,\cdot).
\end{equation*}

For $0\leq k\leq n$, let $\Phi_k(\Omega,\cdot)$ be the hyperbolic
support measure defined by \eqref{eq:tube}, regarded as a measure
on $\mathcal N_\Omega$. On $\mathcal I$, write
$\Phi_k=(J_\Omega)_\#\Phi_k(\Omega,\cdot)$. Since $J_\Omega$
preserves the base point, the projection $\Pi(x,y)=x$ satisfies
$\Pi_\#\Phi_k=V_k(\Omega,\cdot)$. The fiber of $\Pi$ at
$x\in\partial\Omega$ is $\{x\}\times Y_x$, where the exposed
face $Y_x$ of $\Omega^\circ$ is
\begin{equation*}
 Y_x=\{y\in\Omega^\circ:y(x)=1\}
     =\mathcal N_x\cap\{p:h_\Omega(p)=1\}.
\end{equation*}

For $(x,y)\in\mathcal I$, let $(x,\nu)$ be the corresponding pair
in $\mathcal N_\Omega$. The covector $p_{x,\nu}$ is the
$g_x^{-1}$-unit representative of $[y]$. By \eqref{eq:connection}
and $y(x)=1$, we have
$|y|_{g_x^{-1}}^2=(1-|x|^2)(|y|^2-1)$, where $|y|$ is the
Euclidean dual norm. Thus
\begin{equation}\label{var:contact}
 p_{x,\nu}=w(x,y)y,\qquad
 w(x,y)=\frac{1}{|y|_{g_x^{-1}}}
       =\frac{1}{\sqrt{1-|x|^2}\sqrt{|y|^2-1}}.
\end{equation}
The function $w$ is positive and continuous on the compact set
$\mathcal I$.

\begin{pro}\label{prop:localized-variation}
Let $n\geq1$, $0\leq\ell<n$, and $\varepsilon>0$.
Let $\Omega\subset\mathbb H^{n+1}$ be a compact geodesically
convex body with nonempty interior. Suppose that $\Omega$ is a
local minimizer of
\begin{equation*}
 J(K)=W_{\ell+1}(K)-\psi(W_\ell(K))+\varepsilon W_n(K)
\end{equation*}
among compact geodesically convex bodies with nonempty interior
in the Hausdorff topology. Assume that $\psi$ is differentiable
at $W_\ell(\Omega)$ and that $\gamma=\psi'(W_\ell(\Omega))>0$.

Define finite nonnegative Borel measures $\alpha$ and $\beta$
on $\mathcal I$ by
\begin{equation*}
 d\alpha=\frac{\gamma w}{\binom n\ell}\,d\Phi_{n-\ell},
 \qquad
 d\beta=w\left(
 \frac{d\Phi_{n-\ell-1}}{\binom n{\ell+1}}
 +\varepsilon\,d\Phi_0\right).
\end{equation*}
Then $\alpha(\mathcal I)=\beta(\mathcal I)>0$. Moreover, let
$\chi:\mathbb R^{n+1}\to[0,\infty)$ be continuous and convex, and let
$b:(\mathbb R^{n+1})^*\to\mathbb R$ be of the form
\begin{equation*}
 b(y)=\min_{1\leq i\leq N}\bigl(y(v_i)+c_i\bigr),
 \qquad v_i\in\mathbb R^{n+1},\quad c_i\in\mathbb R,\quad N\geq1.
\end{equation*}
If $b\geq0$ on $\Omega^\circ$, then
\begin{equation}\label{var:product-test}
 \int_{\mathcal I}\chi(x)b(y)\,d\beta
 \leq
 \int_{\mathcal I}\chi(x)b(y)\,d\alpha.
\end{equation}
\end{pro}

\begin{proof}
For small $|t|$, the Klein dilation $\varphi_t(x)=(1+t)x$ gives
convex bodies $\varphi_t(\Omega)$ in the Klein ball with
$d_H(\varphi_t(\Omega),\Omega)\to0$. The Klein coordinate velocity
at $t=0$ is $x$, so \eqref{var:contact} gives the normal speed
$p_{x,\nu}(x)=w(x,y)$.
Local minimality and Lemma~\ref{lem:ambient} give
\begin{equation*}
 0=\left.\frac{d}{dt}\right|_{t=0}J(\varphi_t(\Omega))
  =\beta(\mathcal I)-\alpha(\mathcal I).
\end{equation*}
The common mass is positive since
$\beta(\mathcal I)\geq\varepsilon\int_{\mathcal I}w\,d\Phi_0>0$
by \eqref{eq:top}.

Fix $\chi$ and $b$ as in the statement. Choose $C>0$ so large that
$C-\chi(z)c_i\geq0$ for all $z\in\Omega$ and $1\leq i\leq N$,
and $C-\chi(z)b(y)\geq0$ for all $(z,y)\in\Omega\times\Omega^\circ$.
 Define
$F:\Omega\times(\mathbb R^{n+1})^*\to\mathbb R$ by
\begin{equation*}
 F(z,p)=h_\Omega(p)
 \left[C-\chi(z)b\left(\frac{p}{h_\Omega(p)}\right)\right]
 \quad (p\ne0),\qquad F(z,0)=0.
\end{equation*}
For $p\ne0$, the covector $p/h_\Omega(p)$ lies in
$\partial\Omega^\circ$. Thus $F\geq0$, and
$b(p/h_\Omega(p))\geq0$ ensures that $F(\cdot,p)$ is concave.
Since $\chi\geq0$, we have
\begin{equation*}
 F(z,p)=\max_{1\leq i\leq N}
 \bigl\{-\chi(z)p(v_i)+[C-\chi(z)c_i]h_\Omega(p)\bigr\}.
\end{equation*}
The coefficients of $h_\Omega(p)$ in this expression are
nonnegative, so $F(z,\cdot)$ is convex and positively homogeneous
of degree one.

For each $z\in\Omega$, let $L_z$ be the compact convex set
with support function $F(z,\cdot)$ and set
$\Omega_t=\{z+tq:z\in\Omega,\ q\in L_z\}$ for $t\geq0$.
Since $F\geq0$, each $L_z$ contains the origin, so
$\Omega\subset\Omega_t$. Proposition~\ref{prop:variation} gives
geodesic convexity for small $t$, and the uniform bound on the
$L_z$ gives $d_H(\Omega_t,\Omega)\to0$ as $t\downarrow0$.
By \eqref{var:contact} and homogeneity,
$F(x,p_{x,\nu})=w(x,y)F(x,y)$. Using \eqref{var:first} for
$W_{\ell+1}$, $W_\ell$, and $W_n$ in the chain rule for $J$,
we obtain from local minimality
\begin{equation*}
 \begin{aligned}
 0\leq\left.\frac{d}{dt}\right|_{0+}J(\Omega_t)
 &=\int_{\mathcal I}F(x,y)\,d(\beta-\alpha)\\
 &=-\int_{\mathcal I}\chi(x)b(y)\,d(\beta-\alpha).
 \end{aligned}
\end{equation*}
The last equality uses $h_\Omega(y)=1$ on $\mathcal I$ and $\alpha(\mathcal I)=\beta(\mathcal I)$. This proves \eqref{var:product-test}.
\end{proof}
\subsection{An inequality on the polar fibers}\label{subsec:polar-fibre}

Let $\lambda$ and $a_j$ be as in Proposition~\ref{prop:newton}.
Since $\Pi_\#(\Phi_{n-j}/\binom nj)=a_j\lambda$, disintegration
along $\Pi$ with respect to $\lambda$ gives finite measures
$\mu_{j,x}^{\mathrm{pol}}$ on $Y_x$ \cite{ags}*{Theorem~5.3.1}.
For $0\leq j\leq n$ and every bounded Borel function $F$ on $\mathcal I$,
\begin{equation*}
 \frac{1}{\binom nj}\int_{\mathcal I} F\,d\Phi_{n-j}
 =\int_{\partial\Omega}\int_{Y_x} F(x,y)\,
    d\mu_{j,x}^{\mathrm{pol}}(y)\,d\lambda(x).
\end{equation*}
The mass of $\mu_{j,x}^{\mathrm{pol}}$ is $a_j(x)$ for
$\lambda$-almost every $x$. Let $S_x=\{u\in\mathcal N_x:|u|_{g_x^{-1}}=1\}$ be the unit
section of $\mathcal N_x$ with respect to the dual metric $g_x^{-1}$. By \eqref{var:contact}, rescaling
along each conormal ray gives
\begin{equation}\label{fib:unit}
 R_x:Y_x\longrightarrow S_x,\qquad R_x(y)=w(x,y)y.
\end{equation}
We set $\mu_{j,x}=(R_x)_\#\mu_{j,x}^{\mathrm{pol}}$. 
Equivalently, in unit conormal coordinates, disintegrating
$\Phi_{n-j}/\binom nj$ along $(x,u)\mapsto x$ with respect to
$\lambda$ gives the fiber measures $\mu_{j,x}$ on $S_x$.

\begin{corollary}\label{cor:polar-fibre}
Under the assumptions of Proposition~\ref{prop:localized-variation},
set
\begin{equation*}
 d\alpha_x(y)=\gamma w(x,y)\,d\mu_{\ell,x}^{\mathrm{pol}}(y),
 \qquad
 d\beta_x(y)=w(x,y)\bigl(d\mu_{\ell+1,x}^{\mathrm{pol}}(y)
                       +\varepsilon\,d\mu_{n,x}^{\mathrm{pol}}(y)\bigr).
\end{equation*}
There is a $\lambda$-null set $N$ such that, for every exposed
point $x\in\partial\Omega\setminus N$ and every continuous,
nonnegative concave function $f$ on $Y_x$,
\begin{equation}\label{var:polar-fibre}
 \int_{Y_x} f\,d\beta_x
 \leq \int_{Y_x} f\,d\alpha_x .
\end{equation}
\end{corollary}

\begin{proof}
Let $E$ be the set of exposed points of $\Omega$, which is Borel
by \cite{choquetcorsonklee}*{Theorem~1.1}. We first show that if
$\zeta$ and $\xi$ are finite nonnegative Borel measures on $\Omega$ and
$\int\chi\,d\zeta\leq\int\chi\,d\xi$ for every continuous
nonnegative convex function $\chi$ on $\mathbb R^{n+1}$, then
$\zeta|_E\leq\xi|_E$.

Take a nonempty compact set $A\subset E$ and a relatively open
neighborhood $O$ of $A$ in $\Omega$. For each $z\in A$ choose a
linear functional $p_z$ attaining its maximum over $\Omega$ only
at $z$. By compactness of $\Omega\setminus O$, we can choose
$\delta_z>0$ so that
\begin{equation*}
 \chi_z(v)=\max\left\{0,\,1+\frac{p_z(v)-p_z(z)}{\delta_z}\right\}
\end{equation*}
vanishes on $\Omega\setminus O$. On $\Omega$, we have
$0\leq\chi_z\leq1$ and $\chi_z(z)=1$. Fix $0<\eta<1$. By continuity and $\chi_z(z)=1$, the family $\bigl\{\{v:\chi_z(v)>1-\eta\}:z\in A\bigr\}$
is an open cover of $A$. By compactness of $A$, this cover admits a finite subcover. Set $\chi(v)=\max_{1\leq i\leq m}\chi_{z_i}(v)$. Then $\chi$ is
continuous, nonnegative, and convex on $\mathbb R^{n+1}$.
Since $(1-\eta)\mathbf1_A\leq\chi\leq\mathbf1_O$ on $\Omega$,
the assumed integral inequality gives
\begin{equation*}
 (1-\eta)\zeta(A)\leq\int\chi\,d\zeta
 \leq\int\chi\,d\xi\leq\xi(O).
\end{equation*}
Both $\zeta$ and $\xi$ are Radon measures on the compact metric
space $\Omega$. Letting $\eta\downarrow0$ and using the outer
regularity of $\xi$ gives $\zeta(A)\leq\xi(A)$. The inner
regularity of $\zeta$ then gives $\zeta|_E\leq\xi|_E$.

For any concave polyhedral function $b$ admissible in
\eqref{var:product-test}, apply this comparison to
$\zeta=\Pi_\#(b\beta)$ and $\xi=\Pi_\#(b\alpha)$, where $\alpha$
and $\beta$ are as in Proposition~\ref{prop:localized-variation}.
By disintegration over $\lambda$, the measure inequality gives
\begin{equation}\label{var:polyhedral-fibre}
 \int_{Y_x}b\,d\beta_x\leq\int_{Y_x}b\,d\alpha_x
 \qquad\text{for $\lambda$-almost every $x\in E$.}
\end{equation}
Let $\mathcal G_{\mathbb Q}$ be the family of functions $b$ that are
nonnegative on $\Omega^\circ$ and have the form
\[
 b(y)=\min_{1\leq i\leq m}\bigl(y(v_i)+c_i\bigr),
\]
where $m\geq1$, each $v_i$ has rational coordinates in the chosen
Klein chart, and each $c_i$ is rational. Since this family is countable,
we can choose a single $\lambda$-null set $N$, including the exceptional
sets for the disintegrations, such that \eqref{var:polyhedral-fibre} holds for every
$x\in E\setminus N$ and every $b\in\mathcal G_{\mathbb Q}$.

Fix $x\in E\setminus N$ and let $f$ be continuous,
nonnegative, and concave on $Y_x$. Apply \cite{rockafellar}*{Theorem~12.1} to $-f$ extended by $+\infty$ outside $Y_x$. On $Y_x$, $f$ is the pointwise infimum of
its affine majorants defined on the whole covector space. By continuity of $f$ and compactness of $Y_x$, for each
$\eta>0$ finitely many such majorants $\ell_i$ satisfy
$f\leq\min_i\ell_i\leq f+\eta$ on $Y_x$.

The functions $\ell_i+2\eta$ are at least $\eta$ on a relative
neighborhood of $Y_x$ in $\Omega^\circ$. On the complement of this
neighborhood, if nonempty, $1-y(x)$ has a positive minimum.
Thus we can choose $M_i\geq0$ such that
\[
 L_i(y)=\ell_i(y)+2\eta+M_i(1-y(x))\geq\eta
 \qquad(y\in\Omega^\circ).
\]
Choose affine functions $\widetilde L_i$ with rational coefficients
such that $|\widetilde L_i-L_i|\leq\eta/2$ on $\Omega^\circ$.
The bound $L_i\geq\eta$ gives $\widetilde L_i\geq\eta/2$ on
$\Omega^\circ$. Thus $b=\min_i\widetilde L_i$ belongs to
$\mathcal G_{\mathbb Q}$. On $Y_x$, the identity $L_i=\ell_i+2\eta$
and the bounds $f\leq\min_i\ell_i\leq f+\eta$ give $|b-f|\leq4\eta$.
Since the fiber measures are finite, we can let $\eta\downarrow0$
in \eqref{var:polyhedral-fibre} to obtain \eqref{var:polar-fibre}.
\end{proof}
\section{Concavity of the curvature quotient on the conormal cone}
\label{sec:fibre}

We use the measure $\lambda$ and the curvature densities $a_j$ of
Proposition~\ref{prop:newton}. Continuous concave approximations of the
curvature quotient on the polar faces allow us to apply
\eqref{var:polar-fibre} to that quotient. Together with Newton's
inequalities \eqref{newton:projected}, this gives the following estimate.

\begin{pro}\label{prop:total}
Fix $0\leq\ell<n$ and $\varepsilon>0$. Let
$\Omega\subset\mathbb H^{n+1}$ be a compact geodesically convex body
with nonempty interior. Suppose that $\Omega$ locally minimizes
\begin{equation*}
 J(K)=W_{\ell+1}(K)-\psi(W_\ell(K))+\varepsilon W_n(K)
\end{equation*}
among compact geodesically convex bodies in the Hausdorff topology.
Assume that $\psi$ is differentiable at $W_\ell(\Omega)$ and that
$\gamma=\psi'(W_\ell(\Omega))>0$. Then
\begin{equation}\label{fib:target}
 a_n(x)\leq\gamma^{n-\ell}a_\ell(x)
 \qquad\text{for $\lambda$-almost every $x\in\partial\Omega$.}
\end{equation}
In particular, $T_n(\Omega)\leq\gamma^{n-\ell}T_\ell(\Omega)$.
\end{pro}

\begin{remark}\label{rem:lower-order-variation}
For $1\leq\ell<n$, the conclusions of Proposition~\ref{prop:total}
also hold for a local minimizer of $J+\theta W_{\ell-1}$, where
$\theta\geq0$, under the same hypotheses on $\Omega$, $\psi$ and
$\varepsilon$. Indeed, in the proof of
Proposition~\ref{prop:localized-variation} the first variation
replaces $\beta$ by
\begin{equation*}
 \widetilde\beta=\beta+
       \frac{\theta w}{\binom n{\ell-1}}\,\Phi_{n-\ell+1}.
\end{equation*}
The two-sided dilation gives
$\widetilde\beta(\mathcal I)=\alpha(\mathcal I)$.
Using this identity to cancel the constant in the one-sided
variation, we obtain
\begin{equation*}
 \int_{\mathcal I}\chi b\,d\beta
 \leq\int_{\mathcal I}\chi b\,d\widetilde\beta
 \leq\int_{\mathcal I}\chi b\,d\alpha ,
\end{equation*}
where the first inequality uses $\chi b\geq0$.
Thus \eqref{var:product-test} holds for the original $\alpha$ and
$\beta$. The proof of Corollary~\ref{cor:polar-fibre} gives
\eqref{var:polar-fibre} for their conditional measures, and the proof of
Proposition~\ref{prop:total} applies.
\end{remark}

\subsection{Conditional support measures}

Let $\Omega\subset\mathbb H^{n+1}$ be a compact geodesically convex
body with nonempty interior, and choose the Klein origin in
$\operatorname{int}\Omega$. We use the Euclidean support function
$h_\Omega$ and the conormal cone $\mathcal N_x$ in
\eqref{eq:conormal-cone}.

We use the unit conormal sections $S_x$ and the measures
$\mu_{j,x}=(R_x)_\#\mu_{j,x}^{\mathrm{pol}}$ from
Section~\ref{subsec:polar-fibre}.

Let $G$ be an inner product on $\mathbb R^{n+1}$. On the covector space
$(\mathbb R^{n+1})^*$, let $\mathcal H^k_{G^{-1}}$ denote
$k$-dimensional Hausdorff measure induced by the dual norm
$|\cdot|_{G^{-1}}$. For $x\in\partial\Omega$, set
\begin{equation*}
  \sigma_x^G
 =\mathcal H^{\dim\mathcal N_x-1}_{G^{-1}}
   |_{\{u\in\mathcal N_x:\,|u|_{G^{-1}}=1\}}.
\end{equation*}

\begin{lemma}\label{lem:conditional}
For $\lambda$-almost every $x\in\partial\Omega$ with $a_n(x)>0$,
the following assertions hold.
\begin{enumerate}
\item The point $x$ is exposed, and
\begin{equation}\label{fib:dom}
 \mu_{j,x}\ll\mu_{n,x}\sim\sigma_x^{g_x},
 \qquad 0\leq j\leq n.
\end{equation}

\item For $\sigma_x^{g_x}$-almost every $u\in S_x$, the support function
$h_\Omega:(\mathbb R^{n+1})^*\to\mathbb R$ admits a second-order
Alexandrov expansion at $u$. At each such $u$, let $A(u)$ be the matrix of $D^2h_\Omega(u)$ in a $g_x^{-1}$-orthonormal basis of $(\mathbb R^{n+1})^*$. For $\mu_{n,x}$-almost every $u\in S_x$,
\begin{equation}\label{fib:radii}
 \frac{d\mu_{j,x}}{d\mu_{n,x}}(u)
 =\frac{\sigma_{n-j}(A(u))}{\binom nj},
 \qquad 0\leq j\leq n.
\end{equation}
\end{enumerate}
\end{lemma}

\begin{proof}
Let $I$ be the standard Euclidean metric in Klein coordinates.
We first represent each conormal ray by its $I$-unit covector,
identifying vectors and covectors by $I$. The metric comparison
following Lemma~\ref{lem:metric} gives $\Phi_k\sim\Phi_k^I$, hence
$V_k\sim V_k^I$ by projection. Choose a Borel version of $a_n$ and set
$X=\{x\in\partial\Omega:a_n(x)>0\}$. Since $V_0=a_n\lambda$,
we have $V_0^I|_X\sim\lambda|_X$.

Fix $t>0$ and put $T_t=\{z:0<\operatorname{dist}_I(z,\Omega)<t\}$.
Let $P$ be the Euclidean metric projection onto $\Omega$.
For $z\in T_t$, set $u(z)=(z-P(z))/|z-P(z)|$ and
$\tau_t=(P,u)_\#(dz|_{T_t})$. The local Steiner formula
\cite{schneider}*{\S4.2} gives
\begin{equation*}
 \tau_t=\sum_{j=0}^n\frac{t^{j+1}}{j+1}\Phi_{n-j}^I.
\end{equation*}
All coefficients are positive, so $\Pi_\#\tau_t\sim\lambda$ and
$\Phi_{n-j}\ll\tau_t$ for every $j$.

Define the convex function $b$ on $\mathbb R^{n+1}$ by
\begin{equation*}
 b(z)=\max_{v\in\Omega}\left(z\cdot v-\frac12|v|^2\right)
     =\frac12|z|^2-\frac12\min_{v\in\Omega}|z-v|^2 .
\end{equation*}
For each $z\in\mathbb R^{n+1}$, the function
$v\mapsto z\cdot v-\tfrac12|v|^2$ attains its maximum over
$\Omega$ uniquely at $v=P(z)$.
By Danskin's theorem \cite{bonnansshapiro}*{Theorem~4.1},
$\nabla b=P$. For each $x\in\partial\Omega$, a point $z$
satisfies $P(z)=x$ if and only if $z-x\in\mathcal N_x$. Thus
\begin{equation*}
 (\partial b)^{-1}(x)=P^{-1}(x)=x+\mathcal N_x .
\end{equation*}
Since $\Pi(P(z),u(z))=P(z)$ and
$\Pi_\#\Phi_{n-j}^I=V_{n-j}^I$, the local Steiner formula gives
\begin{equation*}
 P_\#(dz|_{T_t})=\Pi_\#\tau_t
 =\sum_{j=0}^n\frac{t^{j+1}}{j+1}V_{n-j}^I.
\end{equation*}
Disintegrating $dz|_{T_t}$ along $P$ and $\tau_t$ along $\Pi$
with respect to the common base measure $\Pi_\#\tau_t$ gives
the fiber measures $\kappa_x$ on $P^{-1}(\{x\})\cap T_t$ and
$\tau_{t,x}$ on $\{u\in\mathcal N_x:|u|=1\}$, respectively. By \cite{ciosmakF}*{Corollary~5.9} applied to $b$ with weight
$\mathbf1_{T_t}$, for $\Pi_\#\tau_t$-almost every $x$,
\begin{equation*}
 \kappa_x\sim\mathcal H^d|_{(x+\mathcal N_x)\cap T_t},
 \qquad d=\dim\mathcal N_x.
\end{equation*}
By uniqueness of disintegration and the volume formula in polar
coordinates $z=x+r\theta$, with volume element
$r^{d-1}\,dr\,d\sigma_x^I(\theta)$ for $0<r<t$, we obtain
\begin{equation*}
 \tau_{t,x}=u_\#\kappa_x
 \sim u_\#\bigl(\mathcal H^d|_{(x+\mathcal N_x)\cap T_t}\bigr)
 =\frac{t^d}{d}\,\sigma_x^I
 \sim\sigma_x^I.
\end{equation*}

Let $du$ be the area measure on the Euclidean unit sphere
$\mathbb S^n\subset(\mathbb R^{n+1})^*$. In our normalization, the
projection $(x,u)\mapsto u$ pushes $\Phi_0^I$ forward to $du$
\cite{schneider}*{\S4.2}. For $du$-almost every $u\in\mathbb S^n$,
$Dh_\Omega(u)$ is the unique support point of $\Omega$ in direction
$u$ \cite{schneider}*{Corollary~1.7.3}. Hence
\begin{equation*}
 \Phi_0^I=(Dh_\Omega,id)_\#du,
 \qquad V_0^I=(Dh_\Omega)_\#du,
\end{equation*}
and $V_0^I$ is concentrated on exposed points.

Let $A=\{p:1<|p|<2\}$, $dm(p)=\mathbf1_A(p)\,dp$, and
$R(p)=p/|p|$. Since $Dh_\Omega(ru)=Dh_\Omega(u)$ for $r>0$,
polar coordinates give
\begin{equation*}
 (Dh_\Omega,R)_\#m=c_A\Phi_0^I,\qquad
 (Dh_\Omega)_\#m=c_AV_0^I,\qquad
 c_A=\int_1^2 r^n\,dr .
\end{equation*}
Let $\rho_x$ be the fiber measures of $m$ under $Dh_\Omega$
over $(Dh_\Omega)_\#m$. Applying
\cite{ciosmakF}*{Corollary~5.9} to $h_\Omega$ on the whole
covector space with $m$ gives
\begin{equation*}
 (\partial h_\Omega)^{-1}(x)=\mathcal N_x,\qquad
 \rho_x\sim\mathcal H^d|_{A\cap\mathcal N_x}
 \quad\text{for $V_0^I$-almost every $x$.}
\end{equation*}
Disintegrating $(Dh_\Omega,R)_\#m=c_A\Phi_0^I$ and using
polar coordinates $p=r\theta$ with $1<r<2$, we obtain
\begin{equation*}
 \Phi_0^I(dx,du)=K_x^I(du)\,dV_0^I(x),\qquad
 K_x^I=R_\#\rho_x\sim
 \left(\int_1^2 r^{d-1}\,dr\right)\sigma_x^I .
\end{equation*}

The Euclidean and hyperbolic unit sections of $\mathcal N_x$
are related by the radial map
\begin{equation*}
 T_x:\{u\in\mathcal N_x:|u|=1\}\longrightarrow S_x,
 \qquad T_x(u)=\frac{u}{|u|_{g_x^{-1}}}.
\end{equation*} 
Applying \cite{ciosmakF}*{Lemma~2.2} to $\Phi_{n-j}\ll\tau_t$ and
$\Phi_0\sim\Phi_0^I$ gives, for $\lambda$-almost every $x\in X$,
\begin{equation*}
 \begin{aligned}
 (T_x^{-1})_\#\mu_{j,x}
 &\ll \frac{d(\Pi_\#\tau_t)}{d\lambda}(x)\,\tau_{t,x}
 \sim\sigma_x^I, \qquad 0\leq j\leq n,\\
 (T_x^{-1})_\#\mu_{n,x}
 &\sim\frac{dV_0^I}{d\lambda}(x)\,K_x^I
 \sim\sigma_x^I .
 \end{aligned}
\end{equation*}
The map $T_x$ is bi-Lipschitz, and pushing forward gives
\begin{equation*}
 \mu_{j,x}\ll\mu_{n,x}
 \sim(T_x)_\#\sigma_x^I\sim\sigma_x^{g_x}.
\end{equation*}
The equivalence $V_0^I|_X\sim\lambda|_X$ shows that $x$ is exposed
for $\lambda$-almost every $x\in X$. This proves \eqref{fib:dom}.

For $\lambda$-almost every $x\in X$, it suffices by \eqref{fib:dom}
to show that $h_\Omega$ admits a second-order expansion at
$\mu_{n,x}$-almost every $u\in S_x$. By Alexandrov's theorem, homogeneity and polar coordinates,
$h_\Omega$ admits a second-order expansion on the Euclidean
unit sphere outside a $du$-null Borel set $N$. Set
\begin{equation*}
 B=\{(x,[p]):p/|p|\in N\},\qquad
 B_x=\{u\in S_x:(x,[u])\in B\}.
\end{equation*}
Since $\Phi_0^I(B)=du(N)=0$ and $\Phi_0\sim\Phi_0^I$,
\begin{equation*}
 0=\Phi_0(B)
   =\int_{\partial\Omega}\mu_{n,x}(B_x)\,d\lambda(x).
\end{equation*}
By \eqref{fib:dom}, $\sigma_x^{g_x}(B_x)=0$ for
$\lambda$-almost every $x\in X$. By positive homogeneity,
the expansion holds at every $u\in S_x\setminus B_x$.

Let $G$ be a constant metric and  $M=\Pi^{-1}(X)$. By \eqref{fib:dom} we have that $\Phi_k|_M\ll\Phi_0|_M$. Then $\Phi_k^G|_M\ll\Phi_0^G|_M$.
Represent each conormal pair $(x,[p])$ by $(x,u)$ with
$u=p/|p|_{G^{-1}}$. Let $\pi(x,u)=u$ be the direction projection
and let $du_G$ denote the area measure on the unit sphere in the
dual metric $G^{-1}$. In our normalization,
$\pi_\#\Phi_k^G=\binom nk S_k^G(\Omega,\cdot)$ and
$\pi_\#\Phi_0^G=du_G$ \cite{schneider}*{\S4.2}.
The linear map $z\mapsto G^{1/2}z$ is an isometry from
$(\mathbb R^{n+1},G)$ to $(\mathbb R^{n+1},I)$.
For a nonzero covector $p$, let $q$ be the covector satisfying
$q(G^{1/2}z)=p(z)$ for all $z\in\mathbb R^{n+1}$.
If $h_\Omega$ admits a second-order expansion at $p$, the chain
rule and homogeneity give
\begin{equation*}
 A_G(p):=D^2h_{G^{1/2}\Omega}(q/|q|)
 =|p|_{G^{-1}}\,G^{1/2}D^2h_\Omega(p)G^{1/2}.
\end{equation*}
The Euclidean area measure formula
\cite{hug}*{Lemma~1 and Theorem~2} gives
\begin{equation*}
 d(\pi_\#\Phi_k^G)_{\mathrm{ac}}(u)=\sigma_k(A_G(u))\,du_G.
\end{equation*}
Choose a Borel set $D_G$ of full $du_G$-measure on which
$h_\Omega$ admits a second-order expansion and
$(\pi_\#\Phi_k^G)_{\mathrm{s}}|_{D_G}=0$ for every $k$, where
$\mathrm{s}$ denotes the part singular with respect to $du_G$.
For $u\in D_G$, the unique support pair is
$s_G(u)=(Dh_\Omega(u),u)$. Thus
\begin{equation*}
 \Phi_k^G|_{\pi^{-1}(D_G)}
 =(s_G)_\#\bigl(\sigma_k(A_G(u))\,du_G|_{D_G}\bigr),
 \qquad 0\leq k\leq n .
\end{equation*}
 Since $\pi_\#\Phi_0^G=du_G$, we have
$\Phi_0^G(\pi^{-1}(D_G^c))=du_G(D_G^c)=0$.
By $\Phi_k^G|_M\ll\Phi_0^G|_M$, it follows that
$\Phi_k^G(M\cap\pi^{-1}(D_G^c))=0$.
Restricting the pushforward formula to $M$ gives
\begin{equation}\label{fib:euclidean-density}
 \frac{d(\Phi_k^G|_M)}{d(\Phi_0^G|_M)}(x,[p])
 =\sigma_k(A_G(p))
 \qquad\text{for $\Phi_0^G|_M$-almost every $(x,[p])$.}
\end{equation}

The local metric comparison \eqref{newton:freeze} transfers
the density formula \eqref{fib:euclidean-density} to the
hyperbolic metric. Fix Borel versions $r_k=d(\Phi_k|_M)/d(\Phi_0|_M)$.
Use the finite covers $U_{m,i}$ from the proof of
Proposition~\ref{prop:newton}, with
$\sup_{U_{m,i}}\|g-G_{m,i}\|\leq\delta_m$ and
$\delta_m\downarrow0$. By \eqref{newton:freeze} at orders $k$
and $0$, together with \eqref{fib:euclidean-density}, we have
over each $U_{m,i}$
\begin{equation*}
 e^{-2c\delta_m}\sigma_k(A_{G_{m,i}}(p))
 \leq r_k(x,[p])
 \leq e^{2c\delta_m}\sigma_k(A_{G_{m,i}}(p))
 \qquad (\Phi_0|_M\text{-a.e.}).
\end{equation*}
Since there are only countably many triples $(m,i,k)$, for
$\Phi_0|_M$-almost every $(x,[p])$ the bounds hold for every
$(m,i,k)$ with $x\in U_{m,i}$. Fix such a pair $(x,[p])$ and,
for each $m$, choose $i(m)$ with $x\in U_{m,i(m)}$. Then $G_{m,i(m)}\to g_x$.
For fixed $p$ at which $h_\Omega$ admits a second-order
expansion, $A_G(p)$ depends continuously on $G$.
Letting $m\to\infty$ gives, for $\Phi_0|_M$-almost every $(x,[p])$,
\begin{equation*}
 r_k(x,[p])=\sigma_k(A_{g_x}(p))=\sigma_k(A(u)),
 \qquad u=\frac{p}{|p|_{g_x^{-1}}},
\end{equation*}
where $A_{g_x}(p)$ represents $D^2h_\Omega(u)$ in a
$g_x^{-1}$-orthonormal basis of the covector space. Taking $k=n-j$ and disintegrating
\begin{equation*}
 \frac{1}{\binom nj}\Phi_{n-j}|_M
 =\frac{r_{n-j}}{\binom nj}\,\Phi_0|_M
\end{equation*}
over $\lambda|_X$, we obtain, for $\lambda$-almost every $x\in X$,
\begin{equation*}
 d\mu_{j,x}(u)
 =\frac{r_{n-j}(x,[u])}{\binom nj}\,d\mu_{n,x}(u)
 =\frac{\sigma_{n-j}(A(u))}{\binom nj}\,d\mu_{n,x}(u)
\end{equation*}
by uniqueness of disintegration. This proves \eqref{fib:radii}.
\end{proof}
\subsection{Concavity on the conormal cone}

Let $\Omega$ be a compact geodesically convex body with nonempty
interior in the Klein model, and let $h_\Omega$ be its Euclidean
support function. Fix $x\in\partial\Omega$ and choose a
$g_x^{-1}$-orthonormal basis of $(\mathbb R^{n+1})^*$.
Throughout this subsection, norms, balls, and orthogonal
complements are taken with respect to $g_x^{-1}$.

Write $U=\operatorname{ri}\mathcal N_x$ for the interior of
$\mathcal N_x$ in its linear span. Let $\mathcal D$ be the set
of covectors $p\in U$ at which $h_\Omega$ admits a second-order
expansion in $(\mathbb R^{n+1})^*$. For $p\in\mathcal D$,
let $A(p)$ be the matrix of $D^2h_\Omega(p)$ in the chosen basis.

\begin{lemma}\label{lem:concave}
Assume that $\mathcal D$ has full measure in $U$ with respect
to Lebesgue measure on $\operatorname{span}\mathcal N_x$.
Then the matrices $A(p)$, $p\in\mathcal D$, have a common image
$V$ and a common kernel $V^\perp$.

Set $d=\dim V$. If $d>0$, the restriction $A(p)|_V:V\to V$
is positive definite for every $p\in\mathcal D$. The inverses
$B(p)=(A(p)|_V)^{-1}$ extend continuously from $\mathcal D$ to $U$
as positive definite operators on $V$. The extension satisfies
\begin{equation}\label{fib:super}
 B(p+q)\succeq B(p)+B(q),\qquad B(tp)=tB(p)
 \qquad(p,q\in U,\ t>0).
\end{equation}

For $0\leq\ell<n$ with $s=n-\ell\leq d$, define $Q$ on
$\mathcal D$ by
\begin{equation}\label{fib:quotient}
 Q(p)=\frac{\binom n\ell}{\binom n{\ell+1}}\,
      \frac{\sigma_{s-1}(A(p))}{\sigma_s(A(p))}.
\end{equation}
Then $Q$ extends continuously to $U$. With $Q=0$ on the
relative boundary of $\mathcal N_x$, the extended function
is nonnegative, concave, and positively homogeneous of degree one on $\mathcal N_x$, and is bounded on $\{p\in\mathcal N_x:|p|=1\}$.
\end{lemma}

\begin{proof}
By \eqref{eq:conormal-cone}, $h_\Omega(p)=p(x)$ on $\mathcal N_x$,
with $Dh_\Omega(p)=x$ for $p\in\mathcal D$. Convexity and homogeneity give
\begin{equation*}
 A(p)\succeq0,\qquad t\mathcal D=\mathcal D,\qquad
 A(tp)=t^{-1}A(p)\quad(t>0).
\end{equation*}
For $p,q\in\mathcal D$, the inclusion $p-\delta q\in\mathcal N_x$
holds for all sufficiently small $\delta>0$ because $p\in U$.
Whenever $\delta>0$ satisfies this inclusion, subadditivity gives
\begin{equation*}
 h_\Omega(p+v)-h_\Omega(p)-v(x)
 \leq h_\Omega(\delta q+v)-h_\Omega(\delta q)-v(x).
\end{equation*}
Fix a covector $w$ and set $v=tw$ in the preceding inequality.
Since $p,\delta q\in\mathcal D$, the second-order expansions at
$p$ and $\delta q$ give
\begin{equation*}
 \frac{t^2}{2}\,w^{\mathsf T}\bigl(A(p)-A(\delta q)\bigr)w
 \leq o(t^2).
\end{equation*}
Dividing by $t^2/2$ and letting $t\to0$ gives
$w^{\mathsf T}A(p)w\leq w^{\mathsf T}A(\delta q)w$.
Since $w$ is arbitrary, this yields the matrix inequality
\begin{equation}\label{fib:compare}
 A(p)\preceq A(\delta q)=\delta^{-1}A(q).
\end{equation}
Thus $\ker A(q)\subset\ker A(p)$. Interchanging $p$ and $q$
gives equality. Since the Hessians are symmetric, their images
are the orthogonal complement $V$ of their common kernel.

Assume $d>0$. For $p,q,p+q\in\mathcal D$, expansion of
$h_\Omega(p+q+t(v+w))\leq h_\Omega(p+tv)+h_\Omega(q+tw)$
at $t=0$ gives
\begin{equation*}
 (v+w)^T A(p+q)(v+w)\leq v^T A(p)v+w^T A(q)w .
\end{equation*}
Fix $z\in V$ and set $w=z-v$ for $v\in V$. The right-hand side is strictly
convex in $v$, and differentiation at its minimum gives
$A(p)v=A(q)w=:\xi$. It follows that $z=(B(p)+B(q))\xi$ and
that the minimum is $v^T\xi+w^T\xi=z^T\xi$. Hence
\begin{equation*}
 z^T A(p+q)z
 \leq \min_{\substack{v,w\in V\\v+w=z}}
      \bigl(v^T A(p)v+w^T A(q)w\bigr)
 =z^T(B(p)+B(q))^{-1}z.
\end{equation*}

Hence we have that $B(p+q)\succeq B(p)+B(q)$  and $B(tp)=tB(p)$ for
$p,q,p+q\in\mathcal D$ and $t>0$.

Fix $p_0\in U$ and a relative ball $B_{2r}(p_0)\subset U$.
For distinct $p,q\in\mathcal D\cap B_r(p_0)$, set $t=|p-q|/r$.
Since $p+(p-q)/t$ and $q+(q-p)/t$ lie in $B_{2r}(p_0)$,
we have $(1+t)p-q,(1+t)q-p\in\mathcal N_x$.
Applying \eqref{fib:compare} with $\delta=(1+t)^{-1}$ in both
directions and inverting on $V$ gives
\begin{equation}\label{fib:B-compare}
 \begin{aligned}
 (1+t)^{-1}B(q)&\preceq B(p)\preceq(1+t)B(q),\\
 \|B(p)-B(q)\|&\leq \frac{|p-q|}{r}\,\|B(q)\|.
 \end{aligned}
\end{equation}
Fix $q_0\in\mathcal D\cap B_r(p_0)$. Since $|p-q_0|<2r$,
\eqref{fib:B-compare} gives
\begin{equation}\label{fib:B-uniform}
 \tfrac13 B(q_0)\preceq B(p)\preceq 3B(q_0)
 \qquad(p\in\mathcal D\cap B_r(p_0)).
\end{equation}
Substituting $\|B(q)\|\leq 3\|B(q_0)\|$ from
\eqref{fib:B-uniform} into \eqref{fib:B-compare} shows that $B$
is Lipschitz on $\mathcal D\cap B_r(p_0)$ with constant
$3\|B(q_0)\|/r$. By density, $B$ extends uniquely to a Lipschitz
map on $B_r(p_0)$. The lower bound in \eqref{fib:B-uniform}
passes to the limit, so the extension is positive definite.
The local extensions agree on overlaps by density and
continuity and thus define $B$ on $U$.

By hypothesis, $U\setminus\mathcal D$ is Lebesgue-null in
$\operatorname{span}\mathcal N_x$. For each $p\in U$,
\begin{equation*}
 \{q\in U:p+q\notin\mathcal D\}
 =U\cap\bigl((U\setminus\mathcal D)-p\bigr)
\end{equation*}
is also null. Hence, by Fubini's theorem, $p,q,p+q\in\mathcal D$
for almost every $(p,q)\in U\times U$. Continuity extends
$B(p+q)\succeq B(p)+B(q)$ to every pair in $U\times U$.
For each $t>0$, the identity $B(tp)=tB(p)$ extends from
$\mathcal D$ to $U$ by density. This proves \eqref{fib:super}.

Let $s=n-\ell\leq d$. For positive definite $C$ on $V$, put
\begin{equation*}
 F(C)=\frac{\binom n\ell}{\binom n{\ell+1}}\,
      \frac{\sigma_{d-s+1}(C)}{\sigma_{d-s}(C)},\qquad \sigma_0=1.
\end{equation*}
The function $F$ is concave, positively homogeneous of degree one,
and nondecreasing in the order of quadratic forms
\cite{quotient}*{Theorem~2.6 and Corollary~5.3}. Set $Q(p)=F(B(p))$
on $U$. By \eqref{fib:super} and the properties of $F$,
\begin{equation*}
 Q(p+q)\geq F(B(p)+B(q))\geq Q(p)+Q(q),\qquad Q(tp)=tQ(p).
\end{equation*}
Thus $Q$ is concave on $U$. For $p\in\mathcal D$, we have $\sigma_s(A(p))>0$ since $A(p)$
has at least $s$ positive eigenvalues. Hence
\begin{equation*}
 \frac{\sigma_{d-s+1}(B(p))}{\sigma_{d-s}(B(p))}
 =\frac{\sigma_{s-1}(A(p))}{\sigma_s(A(p))},
\end{equation*}
which proves \eqref{fib:quotient}.

Set $Q=0$ on $\mathcal N_x\setminus U$. If $p\in U$ and
$q\in\mathcal N_x\setminus U$, then $q+tp\in U$ for $0<t<1$.
Superadditivity of $Q$ on $U$ gives
\begin{equation*}
 Q(p+q)\geq Q((1-t)p)+Q(q+tp)\geq(1-t)Q(p).
\end{equation*}
Letting $t\downarrow0$, there holds $Q(p+q)\geq Q(p)=Q(p)+Q(q)$.
If $p,q\in\mathcal N_x\setminus U$, then
$Q(p+q)\geq0=Q(p)+Q(q)$. The extension is positively
homogeneous and hence concave by superadditivity. It is Borel
because it is continuous on $U$ and zero on its complement in
$\mathcal N_x$.

Finally, for the relative ball $B_{2r}(p_0)\subset U$ chosen above,
$p_0-rp\in U$ whenever $p\in\mathcal N_x$ and $|p|=1$. Hence
\begin{equation*}
 0\leq Q(p)\leq\frac{Q(p_0)-Q(p_0-rp)}{r}
 \leq\frac{Q(p_0)}{r}.
\end{equation*}
\end{proof}

\subsection{Proof of Proposition~\ref{prop:total}}

Choose the Klein origin in $\operatorname{int}\Omega$. On
$\{a_n=0\}$, inequality \eqref{fib:target} follows from
$a_\ell\geq0$. By Lemma~\ref{lem:conditional},
$\lambda$-almost every $x$ with $a_n(x)>0$ is exposed.
Fix $x$ with $a_n(x)>0$ satisfying the conclusions of
Lemma~\ref{lem:conditional} and Corollary~\ref{cor:polar-fibre},
as well as \eqref{newton:projected}. 

By \eqref{var:polar-fibre}, we have
\begin{equation*}
 0<\varepsilon\int_{Y_x}w(x,y)\,d\mu_{n,x}^{\mathrm{pol}}(y)
 \leq\beta_x(Y_x)\leq\alpha_x(Y_x)
 =\gamma\int_{Y_x}w(x,y)\,d\mu_{\ell,x}^{\mathrm{pol}}(y).
\end{equation*}
Here the strict positivity follows from $w>0$ and $a_n(x)>0$.
Hence $a_\ell(x)>0$. By \eqref{newton:projected}, $a_j(x)>0$
for every $\ell\leq j\leq n$.

Fix a $g_x$-orthonormal basis of $\mathbb R^{n+1}$ in Klein
coordinates. At each covector $p$ where $h_\Omega$ admits a
second-order expansion, let $A(p)=D^2h_\Omega(p)$ be its
Hessian matrix in the dual basis. Let $U=\operatorname{ri}\mathcal N_x$
be the interior of the conormal cone in $L=\operatorname{span}\mathcal N_x$.

By Lemma~\ref{lem:conditional}, $h_\Omega$ admits a second-order
expansion in the ambient covector space at
$\sigma_x^{g_x}$-almost every $u$ in the $g_x^{-1}$-unit conormal
section $S_x$. By homogeneity, an expansion at $u$ gives an
expansion at $ru$ for every $r>0$. Polar coordinates in $L$
therefore show that the set $\mathcal D$ of expansion points
in $U$ has full measure in $U$ with respect to Lebesgue measure
on $L$. Lemma~\ref{lem:concave} now gives a common rank $d$
for $A(p)$, $p\in\mathcal D$.

The relative boundary of $\mathcal N_x$ has Lebesgue measure zero
in its span, so polar coordinates give
$\sigma_x^{g_x}(S_x\setminus U)=0$. By \eqref{fib:dom}, we have
$\operatorname{rank}A(u)=d$ for $\mu_{n,x}$-almost every $u\in S_x$.
Taking $j=\ell$ in \eqref{fib:radii} gives
\begin{equation}\label{fib:positive}
 \int_{S_x}\sigma_{n-\ell}(A(u))\,d\mu_{n,x}(u)
 =\binom n\ell\,\mu_{\ell,x}(S_x)
=\binom n\ell\,a_\ell(x)>0.
\end{equation}
If $d<n-\ell$, then $\sigma_{n-\ell}(A(u))=0$ for
$\mu_{n,x}$-almost every $u$, contrary to \eqref{fib:positive}.
Thus $d\geq n-\ell$.

For the function $Q$ in Lemma~\ref{lem:concave},
\eqref{fib:quotient} and \eqref{fib:radii} give
\begin{equation*}
 \frac{d\mu_{\ell+1,x}}{d\mu_{n,x}}(u)
 =\frac{\sigma_{n-\ell-1}(A(u))}{\binom n{\ell+1}}
 =Q(u)\frac{\sigma_{n-\ell}(A(u))}{\binom n\ell}
 =Q(u)\frac{d\mu_{\ell,x}}{d\mu_{n,x}}(u)
\end{equation*}
for $\mu_{n,x}$-almost every $u$. Hence
\begin{equation}\label{fib:density}
 d\mu_{\ell+1,x}=Q\,d\mu_{\ell,x}.
\end{equation}

Since $Y_x$ is the section of $\mathcal N_x$ by $p(x)=1$,
we have $\operatorname{ri}Y_x=Y_x\cap U$. By
Lemma~\ref{lem:concave}, $Q$ is continuous on
$\operatorname{ri}Y_x$, and its assigned value on
$\partial_{\mathrm{rel}}Y_x$ is zero.
If $Y_x$ is a singleton, $Q|_{Y_x}$ is continuous and
\eqref{var:polar-fibre} applies directly. Otherwise, to apply
\eqref{var:polar-fibre}, we approximate $Q|_{Y_x}$ by continuous
concave functions. Let $d_{Y_x}(y)$ be the Euclidean distance
from $y$ to the relative boundary of $Y_x$ in its affine hull
$\operatorname{aff}Y_x$, and set
\begin{equation*}
 Q_k(y)=\min\{Q(y),k\,d_{Y_x}(y)\},\qquad y\in Y_x,\quad k\geq1.
\end{equation*}
Since $Y_x$ is convex, the distance $d_{Y_x}$ to its relative
boundary is concave on $Y_x$. Thus $Q_k$ is nonnegative
and concave. It is continuous on $\operatorname{ri}Y_x$, and
$0\leq Q_k\leq k\,d_{Y_x}$ gives continuity at the relative
boundary. Since $d_{Y_x}>0$ on $\operatorname{ri}Y_x$ and $Q=0$
on $\partial_{\mathrm{rel}}Y_x$, we have $Q_k\uparrow Q$ on $Y_x$.
Applying \eqref{var:polar-fibre} to $Q_k$ and passing to the limit
by monotone convergence gives \eqref{var:polar-fibre} for $Q|_{Y_x}$.

By \eqref{fib:unit} and the homogeneity of $Q$,
we have $Q(R_x(y))=w(x,y)Q(y)$ for every $y\in Y_x$.
Using \eqref{fib:unit} and the definitions of the fiber measures
therefore yields
\begin{equation}\label{fib:quotient-test}
 \int_{S_x} Q\,d\mu_{\ell+1,x}
 +\varepsilon\int_{S_x} Q\,d\mu_{n,x}
 \leq\gamma\int_{S_x} Q\,d\mu_{\ell,x}.
\end{equation}
By the Cauchy--Schwarz inequality, \eqref{fib:density}, and
\eqref{fib:quotient-test},
\begin{align*}
 a_{\ell+1}(x)^2
 &=\left(\int_{S_x}Q\,d\mu_{\ell,x}\right)^2\\
 &\leq a_\ell(x)\int_{S_x}Q^2\,d\mu_{\ell,x}
 =a_\ell(x)\int_{S_x}Q\,d\mu_{\ell+1,x}
 \leq\gamma a_\ell(x)a_{\ell+1}(x).
\end{align*}
It follows that $a_{\ell+1}(x)/a_\ell(x)\leq\gamma$.
The adjacent ratios are nonincreasing by \eqref{newton:projected},
so
\begin{equation*}
 \frac{a_n(x)}{a_\ell(x)}
 =\prod_{j=\ell}^{n-1}\frac{a_{j+1}(x)}{a_j(x)}
 \leq\left(\frac{a_{\ell+1}(x)}{a_\ell(x)}\right)^{n-\ell}
 \leq\gamma^{n-\ell}.
\end{equation*}
Integrating \eqref{fib:target} over $\partial\Omega$ gives
\begin{equation*}
 T_n(\Omega)=\int_{\partial\Omega}a_n\,d\lambda
 \leq\gamma^{n-\ell}\int_{\partial\Omega}a_\ell\,d\lambda
 =\gamma^{n-\ell}T_\ell(\Omega). \qedhere
\end{equation*}

\section{Proof of the main theorem}\label{sec:proof}

\subsection{Restriction formulas for quermassintegrals and total mean curvatures}
\label{sec:dimension}

The following lemma relates the quermassintegrals and total mean
curvatures of a convex body in a totally geodesic subspace to the
corresponding quantities in $\mathbb H^{n+1}$.

\begin{lemma}\label{lem:restrict}
Let $P\subset\mathbb H^{n+1}$ be a totally geodesic subspace of
codimension $a$, where $1\leq a\leq n-1$, and let $\Omega\subset P$
be a compact geodesically convex set with nonempty interior
relative to $P$. Denote by $W'_i(\Omega)$ and $T'_i(\Omega)$
the quermassintegrals and total mean curvatures of $\Omega$ in $P$,
respectively, using the definitions of Section~\ref{sec:prelim}
with $n$ replaced by $n-a$. Then, for $a\leq j\leq n$,
\begin{equation}\label{eq:restrict}
  T_j(\Omega)=b_jT'_{j-a}(\Omega),\qquad
  W_j(\Omega)=b_jd_jW'_{j-a}(\Omega).
\end{equation}
Here $z$ is the last coordinate on $\mathbb S^a$,
$\mathbb S^a_+=\{z\geq0\}$ is the upper hemisphere, and
$d\sigma$ is its area measure. The constants are
\begin{equation*}
  b_j=\frac{\binom{n-a}{j-a}}{\binom nj}
      \int_{\mathbb S^a_+}z^{j-a}\,d\sigma,\qquad
  d_j=\frac{\int_{\mathbb S^a_+}z^{j-a+1}\,d\sigma}
            {\int_{\mathbb S^a_+}z^{j-a}\,d\sigma}.
\end{equation*}
Moreover,
\begin{equation*}
  0<d_a\leq\cdots\leq d_n,\qquad
  b_{j+1}d_{j+1}=b_j\quad(a\leq j<n).
\end{equation*}
\end{lemma}

\begin{proof}
For $a\leq j\leq n$, the restriction formula
\cite{bht}*{Proposition~3.1 and Remark~3.2} at order $n-j$ gives
\begin{equation*}
  T_j(\Omega)
  =\frac{\binom{n-a}{j-a}}{\binom nj}
   \frac{\omega_j}{\omega_{j-a}}\,T'_{j-a}(\Omega)
  =b_jT'_{j-a}(\Omega).
\end{equation*}
The last equality follows by taking $p=j-a$ in the hemisphere integral
\begin{equation*}
  I_p:=\int_{\mathbb S^a_+}z^p\,d\sigma
      =\frac{\omega_{a+p}}{\omega_p},\qquad p=0,1,\ldots.
\end{equation*}

Fix $a\leq j\leq n$. Almost every $j$-plane $L$ meeting $P$
intersects it in a $(j-a)$-plane. Under $L\mapsto L\cap P$, the
invariant measure on the $j$-planes meeting $P$ induces a locally
finite isometry-invariant measure on the $(j-a)$-planes of $P$.
By uniqueness of the invariant measure up to scale, the Crofton
definition gives
\begin{equation*}
  W_j(K)=c_jW'_{j-a}(K)
\end{equation*}
for every compact convex set $K\subset P$, with $c_j$ independent
of $K$.

To determine $c_j$, take a ball $D$ of positive radius in a totally
geodesic $(n+1-j)$-dimensional subspace of $P$. The Cauchy--Crofton
formula~\cite{solanesthesis}*{Theorem~2.2.2}, in the normalization
of our quermassintegrals, gives
\begin{equation*}
  \frac{\binom nj}{\omega_{j+1}}\,W_j(D)
  =\frac{\operatorname{vol}_{n+1-j}(D)}{2\pi}
  =\frac{\binom{n-a}{j-a}}{\omega_{j-a+1}}\,W'_{j-a}(D).
\end{equation*}
Thus
\begin{equation*}
  c_j
  =\frac{\binom{n-a}{j-a}}{\binom nj}
   \frac{\omega_{j+1}}{\omega_{j-a+1}}
  =\frac{\binom{n-a}{j-a}}{\binom nj}\,I_{j-a+1}
  =b_jd_j.
\end{equation*}
Taking $K=\Omega$ proves the formula for $W_j$ in \eqref{eq:restrict}.

For $a\leq j<n$, integration by parts gives
$(j+1)I_{j-a+2}=(j-a+1)I_{j-a}$, hence
\begin{equation*}
  \frac{b_{j+1}d_{j+1}}{b_j}
  =\frac{j+1}{j-a+1}\frac{I_{j-a+2}}{I_{j-a}}=1.
\end{equation*}
Finally, $I_p>0$ and Cauchy--Schwarz gives
$I_{p+1}^2\leq I_pI_{p+2}$, so $0<d_a\leq\cdots\leq d_n$.
\end{proof}
\subsection{Proof of Theorem~\ref{thm:main}}\label{sec:rigidity}

Since each $f_j$ is strictly increasing by \eqref{eq:ball}, it
suffices to prove \eqref{eq:main} for consecutive indices.
Suppose that a consecutive inequality fails. Let $\ell$ be the least index for which there exists a compact
geodesically convex body $\Omega_0\subset\mathbb H^{n+1}$ with
nonempty interior satisfying
\begin{equation*}
  W_{\ell+1}(\Omega_0)<\varphi(W_\ell(\Omega_0)),\qquad
  \varphi=f_{\ell+1}\circ f_\ell^{-1}.
\end{equation*}
By Proposition~\ref{prop:min}, we may choose $M,\varepsilon>0$
so that the functional $J$ defined by \eqref{eq:penalty} attains
a negative minimum among all nonempty compact geodesically convex
subsets of $\mathbb H^{n+1}$. Let $\Omega_*$ be a minimizer. Set
$r=f_\ell^{-1}(W_\ell(\Omega_*))$, and let $\gamma$ be given by
\eqref{eq:gamma}. Then $r>0$, and \eqref{eq:strict-lower} and
\eqref{eq:gamma} give
\begin{equation}\label{eq:min-strict-total}
  T_n(\Omega_*)>\gamma^{n-\ell}T_\ell(\Omega_*).
\end{equation}
If $\Omega_*$ has nonempty interior in $\mathbb H^{n+1}$,
Proposition~\ref{prop:total} gives the opposite inequality.

Suppose now that $\Omega_*$ has empty interior, and let $P$ be its
totally geodesic span, of codimension $a\geq1$. If $\ell=0$, then
$W_0(\Omega_*)=0$. If $1\leq\ell<a$, then
$\ell+\dim P<n+1$, so almost every totally geodesic $\ell$-plane
misses $P$ and the Crofton formula gives $W_\ell(\Omega_*)=0$.
Since $W_\ell(\Omega_*)=f_\ell(r)>0$, we must have
$1\leq a\leq\ell<n$.

To apply Proposition~\ref{prop:total} in $P$, we restrict $J$ to
compact convex bodies $K\subset P$ with nonempty interior relative
to $P$. Substituting \eqref{eq:restrict} into \eqref{eq:penalty}
and using $b_{\ell+1}d_{\ell+1}=b_\ell$, we obtain
\begin{equation*}
  \frac{J(K)}{b_\ell}
  =W'_{\ell+1-a}(K)-\psi\bigl(W'_{\ell-a}(K)\bigr)
   +\varepsilon\frac{b_nd_n}{b_\ell}W'_{n-a}(K),
\end{equation*}
where
\begin{equation*}
  \psi(t)=\frac{M}{b_\ell}
  \bigl(1-\exp[-\varphi(b_\ell d_\ell t)/M]\bigr).
\end{equation*} Since $P$ is the totally
geodesic span of $\Omega_*$, the minimizer has nonempty interior
in $P$. Since $\Omega_*$ minimizes $J$ and $b_\ell>0$, it also minimizes
the restricted functional $K\mapsto J(K)/b_\ell$ over all compact
convex bodies $K\subset P$ with nonempty interior relative to $P$.
Moreover,
\begin{equation*}
  \varepsilon\frac{b_nd_n}{b_\ell}>0,\qquad
  \psi'\bigl(W'_{\ell-a}(\Omega_*)\bigr)=\gamma d_\ell>0.
\end{equation*}
Applying Proposition~\ref{prop:total} in $P$ with index $\ell-a$
and using \eqref{eq:restrict}, we obtain
\begin{align*}
  T_n(\Omega_*)
  &=b_nT'_{n-a}(\Omega_*)\\
  &\leq b_n(\gamma d_\ell)^{n-\ell}T'_{\ell-a}(\Omega_*)\\
  &=\gamma^{n-\ell}T_\ell(\Omega_*)
    \prod_{j=\ell+1}^{n}\frac{d_\ell}{d_j}
  \leq\gamma^{n-\ell}T_\ell(\Omega_*),
\end{align*}
where $b_{j-1}=b_jd_j$ and $d_\ell\leq d_j$ for $j\geq\ell$.
This contradicts \eqref{eq:min-strict-total} and proves
\eqref{eq:main}.

We next consider the equality case of \eqref{eq:main}. Write $r_j(K)=f_j^{-1}(W_j(K))$.
We first show that, for $1\leq m<n$,
\begin{equation}\label{eq:equality-descent}
 r_m(K)=r_{m+1}(K)>0\quad\Longrightarrow\quad r_{m-1}(K)=r_m(K).
\end{equation}
Suppose that some $K_0$ satisfies
$r_{m-1}(K_0)<r_m(K_0)=r_{m+1}(K_0)$. Put
\begin{equation*}
 \varphi=f_{m+1}\circ f_m^{-1},\qquad g=f_{m-1}\circ f_m^{-1},
 \qquad D(K)=W_{m+1}(K)-\varphi(W_m(K)).
\end{equation*}
By \eqref{eq:main}, $D\geq0$ on all nonempty compact convex sets,
and $D(K_0)=0$. Since $W_{m-1}(K_0)<g(W_m(K_0))$, we can first
choose $M$ large and then $\eta>0$ small enough that
\begin{equation*}
 h(t)=M(1-e^{-g(t)/M}),\qquad
 b:=h(W_m(K_0))-W_{m-1}(K_0)-\eta W_n(K_0)>0 .
\end{equation*}
We have $0\leq h\leq\min\{M,g\}$, and \eqref{eq:ball} gives
$g'(f_m(r))=\tanh r$ for $r>0$. Hence $0\leq h'(t)\leq1$ for $t>0$.

For $\delta>0$, let $K_\delta$ minimize
\begin{equation*}
 J_\delta(K)=D(K)+\delta\bigl(W_{m-1}(K)-h(W_m(K))+\eta W_n(K)\bigr)
\end{equation*}
over all nonempty compact convex sets. Such a minimizer exists by
the compactness argument in Proposition~\ref{prop:min}, since
$J_\delta(K)\geq\delta(\eta W_n(K)-M)$. The comparison
$J_\delta(K_\delta)\leq J_\delta(K_0)=-\delta b$, together with
$D(K_\delta)\geq0$ and $h\leq\min\{M,g\}$, gives
\begin{equation}\label{eq:equality-gap}
 0\leq D(K_\delta)\leq\delta M,\qquad
 g(W_m(K_\delta))-W_{m-1}(K_\delta)\geq b,\qquad
 W_n(K_\delta)\leq M/\eta .
\end{equation}
In particular, $r=r_m(K_\delta)>0$. Set $c=\coth r$. The derivative
of $\varphi+\delta h$ at $W_m(K_\delta)$ is
$\gamma=c+\delta h'(W_m(K_\delta))>0$.
We claim that
\begin{equation}\label{eq:equality-upper}
 T_n(K_\delta)\leq\gamma^{n-m}T_m(K_\delta).
\end{equation}
If $K_\delta$ has nonempty interior, this is
Remark~\ref{rem:lower-order-variation} with $\ell=m$,
$\theta=\delta$ and $\varepsilon=\delta\eta$.

If $K_\delta$ has empty interior, let $P$ be its totally geodesic
span and let $a$ be the codimension of $P$ in $\mathbb H^{n+1}$.
Since $W_m(K_\delta)>0$, the Crofton formula gives $1\leq a\leq m$.
As in the proof of the inequality in \eqref{eq:main}, we minimize
$J_\delta/b_m$ over compact convex bodies with nonempty interior
in $P$. The global minimality of $K_\delta$ and $b_m>0$ ensure that
$K_\delta$ is a minimizer in this class.

Write $W'_j(K)$ and $T'_j(K)$ for the $j$th quermassintegral and
total curvature of $K$, respectively, computed in the induced
hyperbolic metric on $P$. Set $\varepsilon=\delta\eta b_nd_n/b_m>0$.
By \eqref{eq:restrict} and $b_{m+1}d_{m+1}=b_m$, we have
\begin{equation*}
\begin{aligned}
 \frac{J_\delta(K)}{b_m}
 &=W'_{m+1-a}(K)
   -\frac{(\varphi+\delta h)\bigl(b_m d_m W'_{m-a}(K)\bigr)}{b_m}\\
 &\quad+\varepsilon W'_{n-a}(K)
   +\frac{\delta}{b_m}W_{m-1}(K).
\end{aligned}
\end{equation*}
Since $b_m d_m W'_{m-a}(K_\delta)=W_m(K_\delta)$,
the chain rule gives
\begin{equation*}
 \left.\frac{d}{dt}
 \frac{(\varphi+\delta h)(b_m d_m t)}{b_m}
 \right|_{t=W'_{m-a}(K_\delta)}
 =d_m(\varphi+\delta h)'(W_m(K_\delta))
 =\gamma d_m.
\end{equation*}

If $a<m$, \eqref{eq:restrict} gives
\begin{equation*}
 \frac{\delta}{b_m}W_{m-1}(K)
 =\frac{\delta b_{m-1}d_{m-1}}{b_m}W'_{m-a-1}(K).
\end{equation*}
If $a=m$, the Crofton formula gives $W_{m-1}(K)=0$ for compact
convex bodies $K\subset P$. Applying
Remark~\ref{rem:lower-order-variation} when $a<m$ and
Proposition~\ref{prop:total} when $a=m$, with index $m-a$, yields
\begin{equation*}
 T'_{n-a}(K_\delta)\leq(\gamma d_m)^{n-m}T'_{m-a}(K_\delta).
\end{equation*}
By \eqref{eq:restrict}, $b_{j-1}=b_jd_j$, and $d_m\leq d_j$ for
$j>m$, we obtain \eqref{eq:equality-upper}
\begin{equation*}
\begin{aligned}
 T_n(K_\delta)
 &=b_n T'_{n-a}(K_\delta)\\
 &\leq \frac{b_n}{b_m}(\gamma d_m)^{n-m}T_m(K_\delta)\\
 &=\gamma^{n-m}T_m(K_\delta)
   \prod_{j=m+1}^{n}\frac{d_m}{d_j}\\
 &\leq\gamma^{n-m}T_m(K_\delta).
\end{aligned}
\end{equation*}

Since $r_{n-1}(K_\delta)\geq r$, \eqref{eq:top} and
\eqref{eq:equality-upper} give
\begin{equation*}
 c^{n-m}T_m(B(r))=T_n(B(r))\leq T_n(K_\delta)
 \leq\gamma^{n-m}T_m(K_\delta).
\end{equation*}
The bound for $W_n$ in \eqref{eq:equality-gap} places all $K_\delta$,
up to isometries, in a fixed ball. By \eqref{eq:rec} and monotonicity
of the quermassintegrals, $T_m(K_\delta)$ is uniformly bounded.
Since $c\geq1$ and $h'\leq1$, we have $\gamma/c\leq1+\delta$, so
\begin{equation*}
 T_m(B(r))-T_m(K_\delta)
 \leq\bigl((1+\delta)^{n-m}-1\bigr)T_m(K_\delta).
\end{equation*}
On the other hand, \eqref{eq:rec} and \eqref{eq:equality-gap} give
\begin{equation*}
\begin{aligned}
 T_m(B(r))-T_m(K_\delta)
 &=m\bigl(g(W_m(K_\delta))-W_{m-1}(K_\delta)\bigr)
   -(n-m)D(K_\delta)\\
 &\geq mb-(n-m)\delta M.
\end{aligned}
\end{equation*}
Letting $\delta\downarrow0$ gives a contradiction and proves
\eqref{eq:equality-descent}.

Now suppose $r_\ell(\Omega)=r_k(\Omega)=r>0$. Applying
\eqref{eq:equality-descent} successively yields
$r_0(\Omega)=r_1(\Omega)=r$. Thus $\Omega$ and $B(r)$ have the same
positive volume and boundary area. Rigidity in the hyperbolic
isoperimetric inequality~\cite{bds} implies that $\Omega$ is a geodesic
ball. Every geodesic ball attains equality by the definition of $f_j$.

Finally, if $W_\ell(\Omega)=0$, equality in \eqref{eq:main} is
equivalent to $W_k(\Omega)=0$. By the Cauchy--Crofton
formula~\cite{solanesthesis}*{Theorem~2.2.2}, $W_k(\Omega)=0$ if and only if $\dim\Omega\leq n-k$.


\begin{bibdiv}
\begin{biblist}
\bibliographystyle{amsplain}
\bib{ags}{book}{
  author={Ambrosio, L.},
  author={Gigli, N.},
  author={Savar\'e, G.},
  title={Gradient Flows in Metric Spaces and in the Space of Probability Measures},
  series={Lectures in Mathematics ETH Z\"urich},
  publisher={Birkh\"auser},
  address={Basel},
  date={2005},
}
\bib{quotient}{article}{
  author={Andrews, B.},
  title={Pinching estimates and motion of hypersurfaces by curvature functions},
  journal={J. Reine Angew. Math.},
  volume={608},
  date={2007},
  pages={17--33},
}

\bib{acw}{article}{
  author={Andrews, B.},
  author={Chen, X.},
  author={Wei, Y.},
  title={Volume preserving flow and Alexandrov--Fenchel type inequalities in hyperbolic space},
  journal={J. Eur. Math. Soc.},
  volume={23},
  date={2021},
  pages={2467--2509},
}

\bib{aw}{article}{
  author={Andrews, B.},
  author={Wei, Y.},
  title={Quermassintegral preserving curvature flow in hyperbolic space},
  journal={Geom. Funct. Anal.},
  volume={28},
  date={2018},
  pages={1183--1208},
}

\bib{awjdg}{article}{
  author={Andrews, B.},
  author={Wei, Y.},
  title={Volume preserving flow by powers of $k$-th mean curvature},
  journal={J. Differential Geom.},
  volume={117},
  date={2021},
  pages={193--222},
}

\bib{bds}{article}{
  author={B\"ogelein, V.},
  author={Duzaar, F.},
  author={Scheven, C.},
  title={A sharp quantitative isoperimetric inequality in hyperbolic $n$-space},
  journal={Calc. Var. Partial Differential Equations},
  volume={54},
  date={2015},
  pages={3967--4017},
}

\bib{bt}{article}{
  author={Bedford, E.},
  author={Taylor, B. A.},
  title={A new capacity for plurisubharmonic functions},
  journal={Acta Math.},
  volume={149},
  date={1982},
  pages={1--40},
}


\bib{valuation}{article}{
  author={Bernig, A.},
  author={Faifman, D.},
  author={Solanes, G.},
  title={Curvature measures of pseudo-Riemannian manifolds},
  journal={J. Reine Angew. Math.},
  volume={788},
  date={2022},
  pages={77--127},
}
\bib{besauwerner}{article}{
  author={Besau, F.},
  author={Werner, E. M.},
  title={The floating body in real space forms},
  journal={J. Differential Geom.},
  volume={110},
  number={2},
  date={2018},
  pages={187--220},
}
\bib{bonnansshapiro}{article}{
  author={Bonnans, J. F.},
  author={Shapiro, A.},
  title={Optimization problems with perturbations: a guided tour},
  journal={SIAM Rev.},
  volume={40},
  number={2},
  date={1998},
  pages={228--264},
}
\bib{bgl}{report}{
  author={Brendle, S.},
  author={Guan, P.},
  author={Li, J.},
  title={An inverse curvature type hypersurface flow in space forms},
  status={preprint},
}

\bib{bht}{article}{
  author={B\"uhler, T.},
  author={Hug, D.},
  author={Th\"ale, C.},
  title={Visibility and intersection density for Boolean models in hyperbolic space},
  journal={Math. Nachr.},
  volume={299},
  date={2026},
  pages={2542--2565},
}

\bib{crscheuer}{article}{
  author={Cabezas-Rivas, E.},
  author={Scheuer, J.},
  title={The quermassintegral-preserving mean curvature flow in the sphere},
  journal={Anal. PDE},
  volume={17},
  date={2024},
  pages={3589--3621},
}

\bib{caravenna}{article}{
  author={Caravenna, L.},
  author={Daneri, S.},
  title={The disintegration of the Lebesgue measure on the faces of a convex function},
  journal={J. Funct. Anal.},
  volume={258},
  date={2010},
  pages={3604--3661},
}

\bib{cgls}{article}{
  author={Chen, C.},
  author={Guan, P.},
  author={Li, J.},
  author={Scheuer, J.},
  title={A fully-nonlinear flow and quermassintegral inequalities in the sphere},
  journal={Pure Appl. Math. Q.},
  volume={18},
  date={2022},
  pages={437--461},
}

\bib{chensun}{article}{
  author={Chen, M.},
  author={Sun, J.},
  title={Alexandrov--Fenchel type inequalities in the sphere},
  journal={Adv. Math.},
  volume={397},
  date={2022},
  pages={108203},
}
\bib{choquetcorsonklee}{article}{
  author={Choquet, G.},
  author={Corson, H.},
  author={Klee, V.},
  title={Exposed points of convex sets},
  journal={Pacific J. Math.},
  volume={17},
  number={1},
  date={1966},
  pages={33--43},
}

\bib{ciosmakF}{report}{
  author={Ciosmak, K. J.},
  title={Measure contraction property on isometric leaves and monotone fibres},
  status={preprint, arXiv:2609.21510v1},
}
\bib{hug}{article}{
  author={Colesanti, A.},
  author={Hug, D.},
  title={Hessian measures of convex functions and applications to area measures},
  journal={J. Lond. Math. Soc. (2)},
  volume={71},
  date={2005},
  pages={221--235},
}
\bib{clm}{article}{
  author={Colesanti, A.},
  author={Ludwig, M.},
  author={Mussnig, F.},
  title={The Hadwiger theorem on convex functions, III: Steiner formulas and mixed Monge--Amp\`ere measures},
  journal={Calc. Var. Partial Differential Equations},
  volume={61},
  number={5},
  date={2022},
  pages={Paper No.~181, 37~pp.},
}
\bib{dinew}{article}{
  author={Dinew, S.},
  author={Lu, C. H.},
  title={Mixed Hessian inequalities and uniqueness in the class
    $\mathcal E(X,\omega,m)$},
  journal={Math. Z.},
  volume={279},
  number={3--4},
  date={2015},
  pages={753--766},
}
\bib{dudagehriglorenz}{article}{
  author={Duda, S.},
  author={Gehrig, E.},
  author={Lorenz, T.},
  title={External ellipsoidal approximations for set evolution equations},
  journal={J. Optim. Theory Appl.},
  volume={192},
  number={3},
  date={2022},
  pages={759--798},
}
\bib{fu1998}{article}{
  author={Fu, J. H. G.},
  title={Some remarks on {L}egendrian rectifiable currents},
  journal={Manuscripta Math.},
  volume={97},
  number={2},
  date={1998},
  pages={175--187},
}
\bib{gallego}{article}{
  author={Gallego, E.},
  author={Solanes, G.},
  title={Integral geometry and geometric inequalities in hyperbolic space},
  journal={Differential Geom. Appl.},
  volume={22},
  date={2005},
  pages={315--325},
}

\bib{gaohugschneider}{article}{
  author={Gao, F.},
  author={Hug, D.},
  author={Schneider, R.},
  title={Intrinsic volumes and polar sets in spherical space},
  journal={Math. Notae},
  volume={41},
  date={2003},
  pages={159--176},
}

\bib{gewangwu}{article}{
  author={Ge, Y.},
  author={Wang, G.},
  author={Wu, J.},
  title={Hyperbolic Alexandrov--Fenchel quermassintegral inequalities II},
  journal={J. Differential Geom.},
  volume={98},
  date={2014},
  pages={237--260},
}
\bib{ghomi2026}{report}{
  author={Ghomi, M.},
  title={Mean curvatures and symmetry of convex hypersurfaces},
  status={arXiv:2608.28410v1},
}

\bib{guanlisurvey}{article}{
  author={Guan, P.},
  author={Li, J.},
  title={Isoperimetric type inequalities and hypersurface flows},
  journal={J. Math. Study},
  volume={54},
  date={2021},
  pages={56--80},
}

\bib{hlodd}{article}{
  author={Hu, Y.},
  author={Li, H.},
  title={Geometric inequalities for hypersurfaces with nonnegative sectional curvature in hyperbolic space},
  journal={Calc. Var. Partial Differential Equations},
  volume={58},
  date={2019},
  pages={Paper No.~55},
}

\bib{huli}{article}{
  author={Hu, Y.},
  author={Li, H.},
  title={Blaschke--Santal\'o type inequalities and quermassintegral inequalities in space forms},
  journal={Adv. Math.},
  volume={413},
  date={2023},
  pages={108826},
}

\bib{hlw}{article}{
  author={Hu, Y.},
  author={Li, H.},
  author={Wei, Y.},
  title={Locally constrained curvature flows and geometric inequalities in hyperbolic space},
  journal={Math. Ann.},
  volume={382},
  date={2022},
  pages={1425--1474},
}

\bib{kohlmann}{article}{
  author={Kohlmann, P.},
  title={Curvature measures and Steiner formulae in space forms},
  journal={Geom. Dedicata},
  volume={40},
  date={1991},
  pages={191--211},
}
\bib{kohlmann1994}{article}{
   author={Kohlmann, P.},
   title={Minkowski integral formulas for compact convex bodies in standard
   space forms},
   journal={Math. Nachr.},
   volume={166},
   date={1994},
   pages={217--228},
}
\bib{lwx}{article}{
  author={Li, H.},
  author={Wei, Y.},
  author={Xiong, C.},
  title={A geometric inequality on hypersurface in hyperbolic space},
  journal={Adv. Math.},
  volume={253},
  date={2014},
  pages={152--162},
}

\bib{luo}{report}{
  author={Luo, T.},
  author={Wei, Y.},
  author={Zhang, R.},
  author={Zhou, R.},
  title={Curvature quotient flows and Alexandrov--Fenchel inequalities in hyperbolic space},
  status={arXiv:2610.35707v1},
}

\bib{lwzsphere}{report}{
  author={Luo, T.},
  author={Wei, Y.},
  author={Zhou, R.},
  title={Alexandrov--Fenchel inequalities for convex domains in the sphere},
  status={arXiv:2609.13829},
}

\bib{makowskischeuer}{article}{
  author={Makowski, M.},
  author={Scheuer, J.},
  title={Rigidity results, inverse curvature flows and Alexandrov--Fenchel type inequalities in the sphere},
  journal={Asian J. Math.},
  volume={20},
  date={2016},
  pages={869--892},
}

\bib{panscheuer}{report}{
  author={Pan, S.},
  author={Scheuer, J.},
  title={The quermassintegral inequalities for horo-convex domains in the sphere},
  status={arXiv:2512.12565},
}

\bib{panyang}{article}{
  author={Pan, S.},
  author={Yang, B.},
  title={A Minkowski type inequality in warped cylinders},
  journal={Trans. Amer. Math. Soc.},
  volume={378},
  date={2025},
  pages={2061--2089},
}
\bib{panasyuk1980}{article}{
  author={Panasyuk, A. I.},
  author={Panasyuk, V. I.},
  title={An equation generated by a differential inclusion},
  journal={Mat. Zametki},
  volume={27},
  number={3},
  date={1980},
  pages={429--437},
  language={Russian},
  translation={
    journal={Math. Notes},
    volume={27},
    number={3},
    date={1980},
    pages={213--218},
  },
}
\bib{rockafellar}{book}{
  author={Rockafellar, R. T.},
  title={Convex Analysis},
  series={Princeton Mathematical Series},
  volume={28},
  publisher={Princeton University Press},
  address={Princeton, NJ},
  date={1970},
}
\bib{santalo}{book}{
  author={Santal\'o, L. A.},
  title={Integral Geometry and Geometric Probability},
  edition={2},
  series={Cambridge Mathematical Library},
  publisher={Cambridge University Press},
  address={Cambridge},
  date={2004},
}

\bib{schneider}{book}{
  author={Schneider, R.},
  title={Convex Bodies: The Brunn--Minkowski Theory},
  edition={2nd expanded ed.},
  publisher={Cambridge University Press},
  date={2014},
}

\bib{solanesthesis}{thesis}{
  author={Solanes, G.},
  title={Integral geometry and curvature integrals in hyperbolic space},
  type={Ph.D. thesis},
  organization={Universitat Aut\`onoma de Barcelona},
  date={2003},
}

\bib{solanes}{article}{
  author={Solanes, G.},
  title={Integral geometry and the Gauss--Bonnet theorem in constant curvature spaces},
  journal={Trans. Amer. Math. Soc.},
  volume={358},
  number={3},
  date={2006},
  pages={1105--1115},
}

\bib{veronelli}{article}{
  author={Veronelli, G.},
  title={Boundary structure of convex sets in the hyperbolic space},
  journal={Monatsh. Math.},
  volume={188},
  date={2019},
  pages={567--586},
}

\bib{wangxia}{article}{
  author={Wang, G.},
  author={Xia, C.},
  title={Isoperimetric type problems and Alexandrov--Fenchel type inequalities in the hyperbolic space},
  journal={Adv. Math.},
  volume={259},
  date={2014},
  pages={532--556},
}

\bib{weixiong}{article}{
  author={Wei, Y.},
  author={Xiong, C.},
  title={Inequalities of Alexandrov--Fenchel type for convex hypersurfaces in hyperbolic space and in the sphere},
  journal={Pacific J. Math.},
  volume={277},
  number={1},
  date={2015},
  pages={219--239},
}

\bib{wyz}{article}{
  author={Wei, Y.},
  author={Yang, B.},
  author={Zhou, T.},
  title={Volume preserving Gauss curvature flow of convex hypersurfaces in the hyperbolic space},
  journal={Trans. Amer. Math. Soc.},
  volume={377},
  date={2024},
  pages={2821--2854},
}

\bib{zhangh}{report}{
  author={Zhang, X.},
  title={Modified support measures and Alexandrov--Fenchel type inequalities for nonsmooth $h$-convex bodies in hyperbolic space},
  status={arXiv:2610.18761v2},
}

\end{biblist}
\end{bibdiv}
\end{document}